# PIML-OFEM: A New Large-Scale Structural Analysis Method Based on Problem-Independent Machine Learning and Overlapping Finite Element Technique

Yilin Guo[1], Chang Liu[1,*], Zongliang Du[1,2], Jin Liu[1], Jingyu Feng[1], Xinyang Zhang[1], Yang Li[1], Tianxing Yang[1], Changyu Shen[1,3,*] and Xu Guo[1,2,*]

[1]*Department of Engineering Mechanics, State Key Laboratory of Structural Analysis, Optimization and CAE Software for Industrial Equipment, Dalian University of Technology, 116023, Dalian, China*

[2]*Dalian Institute of Industrial Software, No. 1 Huoju Road, High-tech Industrial Park, 116085, Dalian, China*

[3]*School of Materials Science and Engineering, The Key Laboratory of Material Processing and Mold of Ministry of Education, Zhengzhou University, 450002, Zhengzhou, China*

**Abstract**

High-resolution structural analysis and rapid design of large-scale heterogeneous structures require both accurate reduced-order models and efficient online computation. Under the traditional multiscale finite element framework, constructing multiscale shape functions online for different material distributions incurs prohibitive computational costs. Meanwhile, in previously developed substructure-based problem-independent machine learning (PIML) techniques, the analysis accuracy is often constrained by the assumed form of displacement distribution on the substructure boundaries. To address these challenges, we propose PIML-OFEM—an overlapping finite element method accelerated by problem-independent machine learning. This method retains only the corner-node degrees of freedom of each substructure and constructs numerical basis functions associated with these degrees of freedom by solving local boundary value problems on extended domains. This reduced-order modeling approach eliminates the need to prescribe displacement assumptions on the target substructure boundaries, significantly enhancing the ability of the numerical basis functions to reproduce local deformations within the substructure. A partition-of-unity overlapping finite element formulation is further introduced to interpolate the independently constructed numerical basis functions into globally continuous basis functions. To avoid solving for the numerical basis

[*] Corresponding authors. E-mail:
c.liu@dlut.edu.cn (Chang Liu), shency@zzu.edu.cn (Changyu Shen), guoxu@dlut.edu.cn (Xu Guo)

functions online, we train a U-Net to learn the deterministic mapping from the Young's modulus distribution within a substructure to its numerical basis functions. Since this learning process does not depend on specific loading conditions or boundary conditions, the learned numerical basis functions can be applied to the solution of any similar structural analysis and optimization problem. Numerical examples demonstrate that the displacements and elemental strain energies obtained by PIML-OFEM are in close agreement with fine-scale finite element results. Its online computational cost is significantly lower than that of direct finite element analysis, while its computational accuracy substantially outperforms the PIML method based on linear boundary interpolation. When embedded in a topology optimization framework, PIML-OFEM enables high-resolution optimization with small filter radii and preserves fine-scale structural features, including local patterns resembling rank-2 microstructures. By integrating mechanics-principled local independent reduction, globally continuous coupling, and problem-independent machine learning techniques, PIML-OFEM establishes an efficient AI-empowered computational paradigm for large-scale heterogeneous structural analysis and high-resolution topology optimization.

## 1 Introduction

Efficient and scalable numerical simulation methods are fundamental to the analysis and design optimization of complex engineering structures. For structures with highly heterogeneous material distributions, complex geometric boundaries, or strong local responses, classical finite element methods typically require extremely high-resolution fine-scale meshes to accurately capture displacement, strain, and strain energy distributions. Particularly in structural topology optimization [1], finite element analysis must be performed repeatedly within the optimization iteration loop, and the computational cost grows rapidly with the number of design degrees of freedom and mesh resolution. Therefore, how to significantly reduce the online solution cost while maintaining fine-scale analysis accuracy is a key challenge in large-scale heterogeneous structural analysis and high-resolution topology optimization.

In recent years, machine learning methods have provided new approaches for accelerating scientific computing and have been applied to problems such as turbulence modeling [2-4], materials design [5-7], and numerical simulation [8-10]. However, many machine learning-enhanced numerical methods still exhibit notable problem dependence. For example, some end-to-end models directly learn the mapping from specific geometries, boundary conditions, and loading forms to the overall physical fields [11-18]; other methods treat neural networks as trainable function approximators and re-optimize network parameters for a given problem, such as the HiDeNN method [19-22] and related works that introduce piecewise polynomials and weak-form residuals into neural representations [23-24]. In topology optimization, machine learning has also been used to directly predict sensitivities, eliminating the need to solve the structural displacement field and thereby reducing computational cost [25-26]. Such methods can achieve high efficiency for specific test cases, but when the structural dimensions, material topology, boundary conditions, or loading forms change, retraining or re-optimization is often required, limiting their generalization capability and model reusability.

To improve the transferability of machine learning models in mechanical analysis, the problem-independent machine learning model has been adopted for large-scale mechanical analysis and topology optimization [27-30]. The core idea of this class of methods is not to learn the global solution of a specific boundary value problem, but rather to learn local numerical operators or multiscale numerical basis functions uniquely determined by the local material distribution. Since these local mappings are decoupled from the overall structural geometry, external loads, and boundary conditions, the trained model can, in principle, be reused for different structural analysis tasks governed by the same class of partial differential equations. Based on this idea, Guo et al. proposed Problem Independent Machine Learning (PIML) paradigm which can significantly reduce the online degrees of freedom through substructure condensation [31-32] and use offline-trained machine learning models to rapidly predict local numerical basis functions, thereby improving the computational efficiency of large-scale structural analysis and topology optimization.

Although PIML methods have been applied to ultra-high-resolution topology optimization [33], lattice structures [34], shell-infill composite structures [35], and engineering structural optimization with complex design domains [36], existing substructure-condensation-based PIML analysis frameworks still face a fundamental accuracy-efficiency trade-off. If each two-dimensional substructure retains only four corner nodes, totaling eight degrees of freedom, it is typically necessary to assume that the substructure boundary displacements can be obtained by linear interpolation of the corner node displacements; this boundary displacement assumption weakens the method's ability to represent complex local deformations, disturbances near voids, and sharp responses near loading points. If accuracy is improved by increasing boundary control degrees of freedom or adopting higher-order boundary interpolation [37], this increases the size of the condensed system, reduces the sparsity of the global stiffness matrix, and expands the output dimension of the machine learning model, thus undermining the original efficiency advantage of the substructure reduction approach. Therefore, how to improve the displacement recovery capability within and near substructure boundaries without increasing boundary degrees of freedom is a core problem that existing PIML substructure methods need to address.

Oversampling techniques provide an effective approach for mitigating the aforementioned boundary layer errors [38-39]. By solving local boundary value problems in an extended region covering the target substructure, one can construct multiscale numerical basis functions adapted to the local heterogeneous material distribution, thereby avoiding the imposition of linear or higher-order interpolation assumptions on the boundary displacements of the target substructure region. However, introducing oversampling shape functions alone is insufficient to form a globally usable reduced-order analysis framework. Because the oversampling shape functions of different substructures are determined by their respective local material distributions, adjacent substructures, while maintaining consistency at shared corner nodes, generally yield different displacement distributions along shared edges, leading to a discontinuous global displacement field and compromising the reliability of derived quantities such as strain, strain energy, and topology optimization sensitivities. Furthermore, if oversampling local problems are solved individually for each substructure during the online stage, the computational cost remains high and is difficult to meet the efficiency demands of large-scale analysis and multiple rounds of optimization iterations.

To address the above issues, this paper proposes a problem-independent machine learning-based overlapping finite element substructure analysis framework, termed PIML-OFEM. This framework first constructs oversampling-based multiscale numerical basis functions for each substructure, retaining only four corner nodes with eight degrees of freedom in two-dimensional problems, so as to enhance the expressive capacity of the low-dimensional reduced space for complex local displacement patterns. Subsequently, this paper adopts the overlapping finite element and partition of unity ideas [40-48] to coordinate the local interpolations of adjacent substructures into a globally continuous displacement field through weighted combination, thereby avoiding displacement jumps on shared edges caused by independently constructed substructure shape functions. Building on this, the paper further follows the PIML framework [27-29] to construct a U-Net prediction model that learns the deterministic mapping from the oversampling region material distribution to the local numerical basis functions, replacing expensive online local boundary value problem solutions with batched neural network inference.

The main contributions of this paper are embodied in the following four aspects. First, we construct oversampled numerical basis functions that do not require imposing linear assumptions on the target substructure boundary displacements, enhancing the recovery capability for complex local deformations while maintaining only four corner nodes with eight degrees of freedom. Second, we introduce an overlapping finite element partition-of-unity compatibility mechanism that embeds the local reduced-order displacement fields of independent substructures into a globally continuous displacement space, enabling the resulting solution to be used for post-processing tasks sensitive to continuity, such as strain, strain energy, and topology optimization sensitivity computation. Third, we establish a U-Net prediction model oriented toward local material distributions, decoupling the machine learning module from overall structural dimensions, boundary conditions, and loading forms, thus demonstrating problem-independent characteristics. Fourth, by exploiting regular Cartesian grids and a finite number of overlapping fine-element types, the integration and local stiffness computations associated with weight functions are preprocessed offline, thereby maintaining the efficiency of the online assembly and solution workflow. Numerical examples demonstrate that the proposed method can achieve analysis accuracy approaching that of fine-scale finite elements across various heterogeneous structures and high-resolution topology optimization problems, while significantly reducing the online solution cost.

The remainder of this paper is organized as follows. Section 2 introduces the substructure reduced-order analysis method based on oversampled numerical basis functions, including classical substructure condensation and the construction of oversampling shape functions without target boundary displacement assumptions. Section 3 describes the global displacement compatibility method based on overlapping finite elements and partition of unity, along with its coupling with the oversampling reduced-order mapping. Section 4 presents the U-Net model for rapid prediction of oversampled numerical basis functions , training data, and loss function. Section 5 validates the accuracy, efficiency, interface compatibility, and problem-independent generalization capability of the proposed method through cantilever beam, complex topology structure, and topology optimization examples. Section 6 concludes the paper and discusses future research directions.

# 2 Substructure Reduced-Order Modeling Based on Oversampled Numerical basis functions

This section first reviews the basic formulation of classical substructure condensation, and then constructs oversampled numerical basis functions for individual substructures on this basis. It should be noted that the focus of this section is to establish the local reduced-order space and its condensed stiffness expression; the global displacement compatibility problem arising from stitching shape functions of multiple independent substructures will be further addressed in Section 3 using overlapping finite elements and the partition of unity method.

## 2.1 Classical Substructure Static Condensation Method

The substructure method can be regarded as a static condensation method. Its core idea is to eliminate the internal degrees of freedom within each subdomain, retaining only the boundary degrees of freedom to participate in the global solution; after the boundary degrees of freedom are obtained, the internal displacements of the substructure are then recovered through numerical basis functions. If the complete substructure boundary degrees of freedom are retained, this condensation process introduces no additional approximation at the discrete equation level, and thus can preserve a solution consistent with the original fine-scale finite element system while reducing the scale of the global solution.

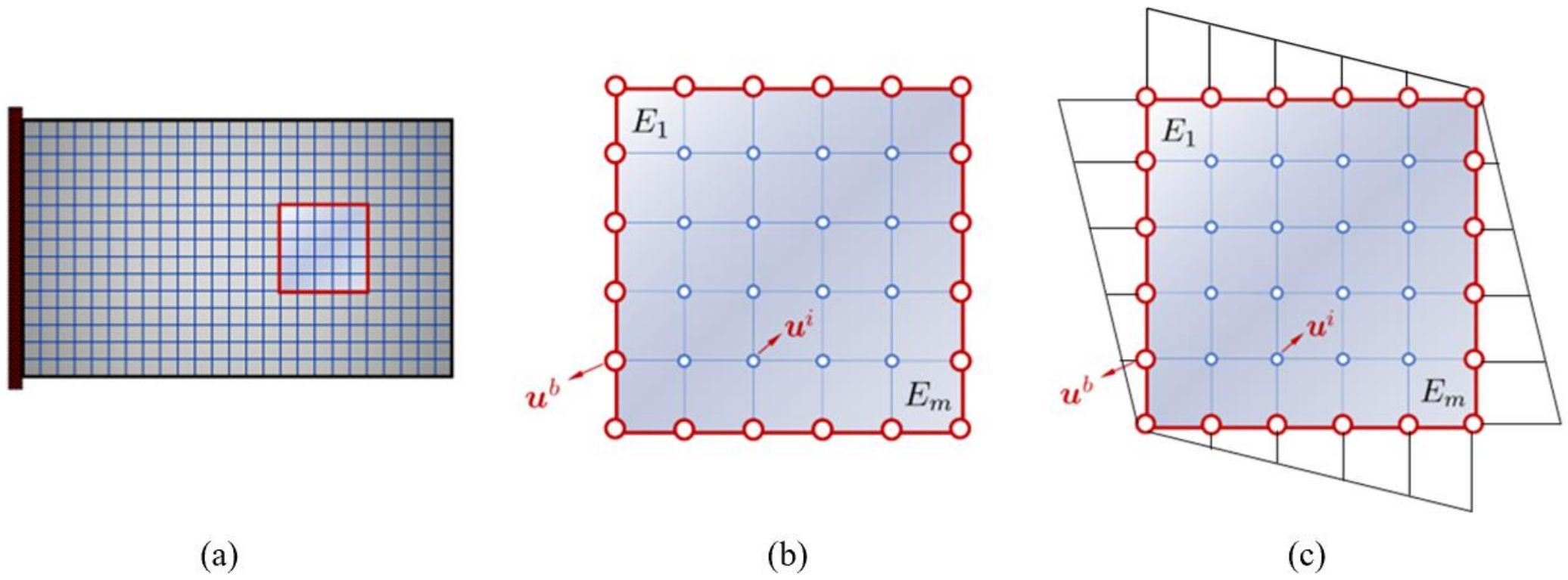

Fig.1. Classical substructure condensation. (a) Decomposition of the analysis domain into substructures. (b) Partition of the nodal degrees of freedom into internal and boundary sets. (c) Linear interpolation of the substructure-boundary displacements from the corner-node degrees of freedom.

Next, we derive the numerical implementation process of the substructure method in the context of linear elasticity, small-deformation finite element analysis. Combining domain decomposition and finite elements, as shown in Fig. 1(a), the static equilibrium equation within each substructure $\Omega^j$ , $j = 1, \ldots, N$ (note that $\Omega^j$ contains multiple fine-scale elements) can be independently expressed as:

$$\mathbf{K}^j \boldsymbol{u}^j = \boldsymbol{f}^j, \tag{1}$$

where $\mathbf{K}^j$ is the stiffness matrix of substructure $\Omega^j$, and $\boldsymbol{u}^j$ and $\boldsymbol{f}^j$ are the nodal displacement vector and load vector of the substructure, respectively. Next, as shown in Fig. 1(b), the nodes of each substructure are divided into interior nodes and boundary nodes, denoted by subscript i and b respectively. Then Eq. (1) can be decomposed as:

$$\mathbf{K}^j \boldsymbol{u}^j = \begin{pmatrix} \mathbf{K}_{\mathrm{bb}}^j & (\mathbf{K}_{\mathrm{ib}}^j)^{\mathsf{T}} \\ \mathbf{K}_{\mathrm{ib}}^j & \mathbf{K}_{\mathrm{ii}}^j \end{pmatrix} \begin{pmatrix} \boldsymbol{u}_{\mathrm{b}}^j \\ \boldsymbol{u}_{\mathrm{i}}^j \end{pmatrix} = \begin{pmatrix} \boldsymbol{f}_{\mathrm{b}}^j \\ \boldsymbol{f}_{\mathrm{i}}^j \end{pmatrix} \tag{2}$$

Without loss of generality, assume that the external loads on the internal nodes of the substructure are zero, i.e. $\boldsymbol{f}_{\mathrm{i}}^j = 0$. This framework remains applicable to cases with mechanical or internal nodal loads; see Reference [37] for details. Substituting into Eq. (2) yields:

$$\boldsymbol{u}_{\mathrm{i}}^j = -(\mathbf{K}_{\mathrm{ii}}^j)^{-1} \mathbf{K}_{\mathrm{ib}}^j \boldsymbol{u}_{\mathrm{b}}^j. \tag{3}$$

Further substituting Eq. (3) into Eq. (2), the condensed stiffness matrix $\mathbf{K}_{\mathrm{s}}^j$ corresponding to the substructure boundary displacements $\boldsymbol{u}_{\mathrm{b}}^{j}$ can be obtained as:

$$\mathbf{K}_{\mathrm{s}}^{j}=\mathbf{K}_{\mathrm{bb}}^{j}-\left(\mathbf{K}_{\mathrm{ib}}^{j}\right)^{\top}\left(\mathbf{K}_{\mathrm{ii}}^{j}\right)^{-1}\mathbf{K}_{\mathrm{ib}}^{j}, \tag{4}$$

the substructure condensed stiffness matrices are then assembled to obtain the global condensed stiffness matrix $\mathbf{K}_{\mathrm{s}}=\sum\left(\mathbf{G}^{j}\right)^{\top}\mathbf{K}_{\mathrm{s}}^{j}\mathbf{G}^{j}$, where $\mathbf{G}^{j}$ is the assembly position matrix of the $j$-th substructure. Further, by solving:

$$\mathbf{K}_{\mathrm{s}}\boldsymbol{u}_{\mathrm{b}}=\boldsymbol{f}_{\mathrm{b}}, \tag{5}$$

the displacement vector $\boldsymbol{u}_{\mathrm{b}}$ of all substructure boundary nodes can be obtained (here $\boldsymbol{f}_{\mathrm{b}}=\sum\mathbf{G}^{j}\boldsymbol{f}_{\mathrm{b}}^{j}$ is the external load vector assembled from $\boldsymbol{f}_{\mathrm{b}}^{j}, j=1,\dots,N$, corresponding to all substructure boundary nodes). Once $\boldsymbol{u}_{\mathrm{b}}$ is obtained, the displacements of the interior nodes of each substructure can be recovered from Equation (3). It is worth noting that the above static condensation process significantly reduces the scale of the linear system of equations that needs to be solved, thereby improving the solution efficiency for large-scale problems.

From Eq. (3), it can be seen that there exists a linear mapping between the internal node displacements and the boundary node displacements. Therefore, this relationship can be written in the form of numerical basis functions:

$$\boldsymbol{u}_{\mathrm{i}}^{j}=\mathbf{N}_{\mathrm{s}}^{j}\boldsymbol{u}_{\mathrm{b}}^{j}, \tag{6}$$

where $\mathbf{N}_{s}^{j}\in\mathbb{R}^{ni\times nb}$ is the multiscale numerical basis function matrix relating the substructure boundary node displacements to the substructure internal node displacements (with $ni$ and $nb$ denoting the total number of internal node degrees of freedom and boundary node degrees of freedom of the substructure, respectively). By comparing Eq. (3) and Eq. (6), it is evident that $\mathbf{N}_{\mathrm{s}}^{j}=-\left(\mathbf{K}_{\mathrm{ii}}^{j}\right)^{-1}\mathbf{K}_{\mathrm{ib}}^{j}$.

After merging the internal and boundary degrees of freedom, the nodal displacements of the entire substructure can be further expressed as a linear function of the boundary degrees of freedom:

$$\boldsymbol{u}^{j}=\left[\left(\boldsymbol{u}_{\mathrm{i}}^{j}\right)^{\top},\left(\boldsymbol{u}_{\mathrm{b}}^{j}\right)^{\top}\right]^{\top}=\left[\mathbf{N}_{\mathrm{s}}^{j\top},\boldsymbol{I}\right]^{\top}\boldsymbol{u}_{\mathrm{b}}^{j}=\mathbf{N}^{j}\boldsymbol{u}_{\mathrm{b}}^{j}, \tag{7}$$

where $\boldsymbol{I}\in\mathbb{R}^{nb\times nb}$ is the identity matrix and $\mathbf{N}^{j}$ maps the retained boundary degrees of freedom to all nodal degrees of freedom in the substructure.

Furthermore, the substructure condensed stiffness matrix $\mathbf{K}_{\mathrm{s}}^{j}$ can be expressed using the numerical basis function matrix as [29]:

$$\mathbf{K}_{\mathrm{s}}^{j} = \left(\mathbf{N}^{j}\right)^{\mathsf{T}} \mathbf{K}^{j} \mathbf{N}^{j}. \tag{8}$$

In large-scale analysis, if all boundary node degrees of freedom are retained, the scale of the condensed system may still be relatively large. To further improve efficiency, a common practice is to retain only the substructure vertex degrees of freedom and assume that the boundary node displacements can be obtained by linear interpolation of the vertex displacements. In this case, analogous to Eq. (7), the condensed stiffness matrix based on the linear boundary-displacement assumption can be written as:

$$\widetilde{\mathbf{K}}_{\mathrm{s}}^{j} = \left(\mathbf{N}^{j}\mathbf{L}\right)^{\mathsf{T}} \mathbf{K}^{j} \left(\mathbf{N}^{j}\mathbf{L}\right) \triangleq \left(\widetilde{\mathbf{N}}^{j}\right)^{\mathsf{T}} \mathbf{K}^{j} \widetilde{\mathbf{N}}^{j}, \tag{9}$$

where $\widetilde{\mathbf{N}}^{j} = \mathbf{N}^{j}\mathbf{L}$ is the multiscale numerical basis function interpolating the substructure vertex displacements to obtain the full nodal displacement vector of the substructure. In this way, the dimension and computational cost of the condensed stiffness matrix can be further reduced. It is also worth noting that this substructure analysis based on the boundary deformation linear assumption can be viewed as performing coarse-mesh (substructure) analysis of the entire structure using first-order finite elements, where the element stiffness matrix is determined by the material distribution within the fine mesh inside the coarse mesh.

## 2.2 Construction of Oversampled Numerical basis functions and Boundary Treatment

The previous section showed that when the substructure method further retains only corner node degrees of freedom, the recovery of internal displacements of the substructure typically relies on the linear boundary-displacement assumption. For PIML substructure methods that use machine learning to predict numerical basis functions, the computational error mainly arises from two sources: the prediction error of the machine learning model for the numerical basis functions, and the boundary representation error caused by the linear boundary-displacement assumption. The former can be mitigated by increasing training samples, enhancing the expressive capacity of the neural network, or introducing physical constraints; the latter stems from the inability of low-order boundary interpolation to capture oscillatory boundary responses induced by local material heterogeneity. Therefore, even if the numerical basis function prediction is sufficiently accurate, the displacement pattern at the boundary of the target substructure may still be constrained by the linear boundary assumption. To alleviate this problem, a natural approach is to introduce oversampling techniques: instead of directly imposing linear assumptions on the boundary displacements of the target substructure, local numerical basis functions are constructed within an extended region covering the target substructure and then restricted to the target substructure. In this way, the boundary displacements of the target substructure can be naturally induced by the fine-scale material distribution and local mechanical response within the extended region, thereby enhancing the expressive capacity of the low-degree-of-freedom substructure space for complex local deformations.

Another way to improve the boundary representation capability is to adopt higher-order boundary displacement assumptions, interpolating more complex substructure boundary deformations by retaining more boundary control degrees of freedom [37]. However, this strategy directly increases the degrees of freedom of the condensed system and may reduce the sparsity of the global stiffness matrix, thereby increasing the cost of solving the linear system. When the substructure size is further increased and more fine-scale degrees of freedom need to be condensed, even higher-order boundary interpolation schemes are often required to describe more complex boundary deformations, which limits the upper bound of overall efficiency improvement. At the same time, higher-order interpolation schemes also expand the scale of numerical basis functions that the machine learning model needs to predict, thereby increasing the training difficulty and potentially affecting prediction accuracy. Based on the above considerations, this paper aims to improve the recovery accuracy of the displacement field within the target substructure by means of oversampling techniques while still retaining only four corner nodes with eight degrees of freedom [38-39]. It should be noted that existing oversampling methods typically require first constructing temporary numerical basis functions and then obtaining the numerical basis function matrix on the target substructure through additional local finite element analysis; moreover, when constructing oscillatory boundary conditions, the consideration of the Poisson effect and the coupling relationship between different displacement components is still insufficient, which may affect the computational accuracy of derived quantities such as strain and stress. This paper will construct complete oversampled numerical basis functions on this basis to improve the local displacement recovery capability while maintaining a low-degree-of-freedom condensed system.

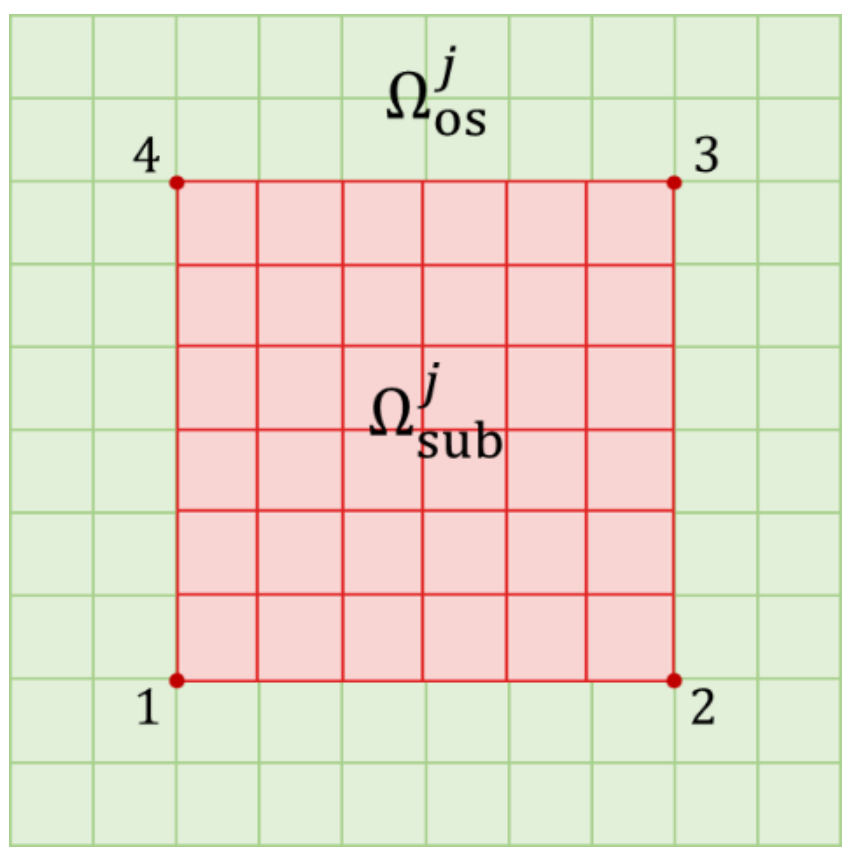

Fig.2. Target substructure $\Omega_{\mathrm{sub}}^{j}$ and its associated oversampling domain $\Omega_{\mathrm{os}}^{j}$.

For each square substructure $\Omega_{\mathrm{sub}}^{j}$, define a larger square oversampling domain $\Omega_{\mathrm{os}}^{j}$ such that $\Omega_{\mathrm{sub}}^{j} \subset \Omega_{\mathrm{os}}^{j}$. Both domains are discretized with square fine-scale elements. In Fig. 2, the red region denotes the target substructure $\Omega_{\mathrm{sub}}^{j}$ and the surrounding green region denotes $\Omega_{\mathrm{os}}^{j} \backslash \Omega_{\mathrm{sub}}^{j}$. The target substructure contains m × m fine-scale elements ($m = 5$ in Fig. 2), and the oversampling domain extends by $l$ element layers on each side ($l = 2$ in Fig. 2).

The objective is to represent all nodal displacements in $\Omega_{\mathrm{sub}}^{j}$ using only the eight displacement degrees of freedom at its four corner nodes, without prescribing a linear or higher-order displacement interpolation along $\partial\Omega_{\mathrm{sub}}^{j}$. To construct the required eight-dimensional local space, auxiliary boundary-value problems are solved over $\Omega_{\mathrm{os}}^{j}$. The resulting displacement basis vectors incorporate the effect of material heterogeneity in the oversampling region and are subsequently restricted to $\Omega_{\mathrm{sub}}^{j}$.

Specifically, the local boundary value problem defined on the oversampling domain $\Omega_{\mathrm{os}}^{j}$ is employed to generate displacement basis vectors exhibiting oscillatory boundary characteristics, which are then restricted to the interior of the target substructure $\Omega_{\mathrm{sub}}^{j}$. The basis vectors obtained in this manner can incorporate the influence of material heterogeneity within the oversampling domain on the boundary response of the target substructure. The construction procedure is described as follows.

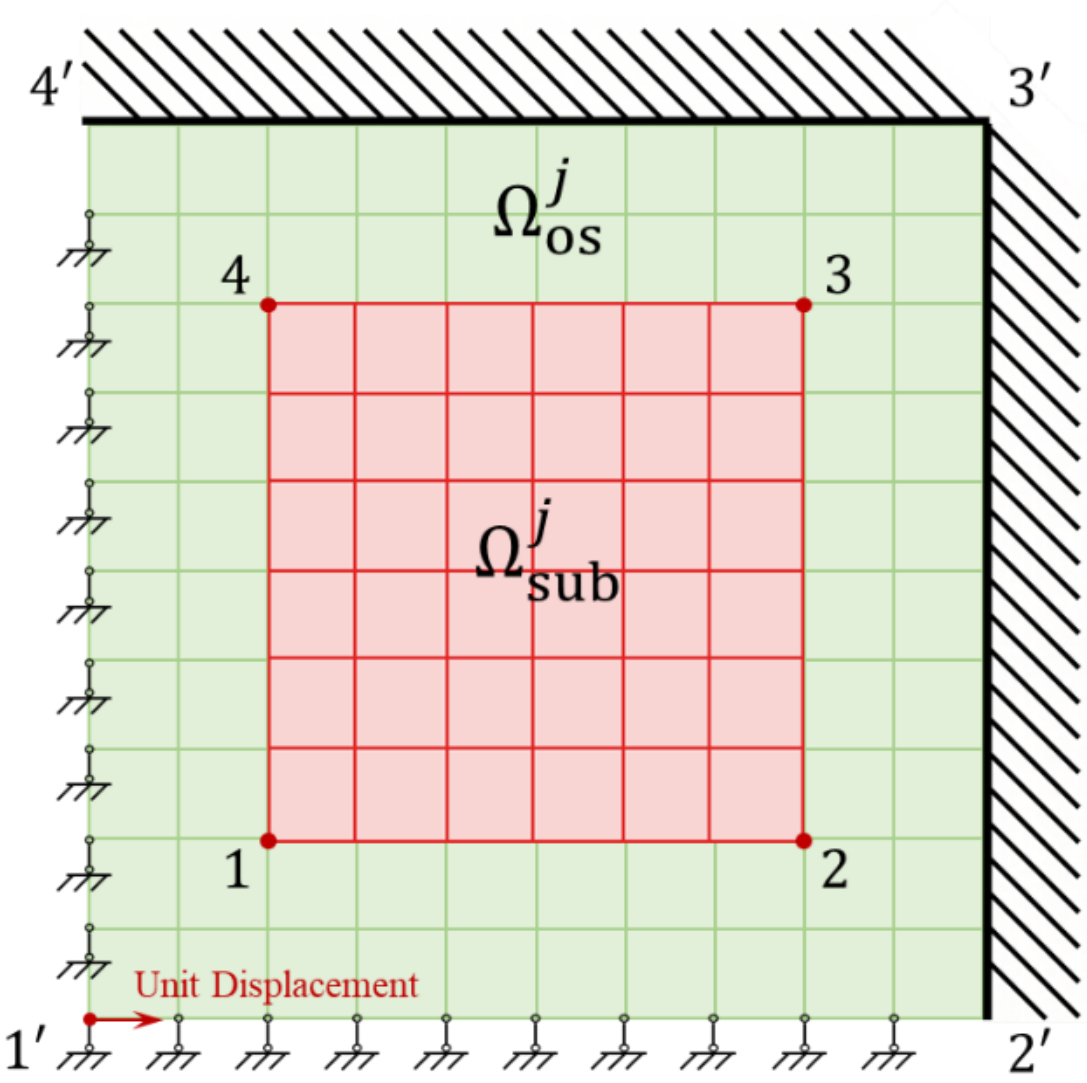

Fig.3. Boundary conditions imposed on the oversampling domain for obtaining $\boldsymbol{\psi}_x^{1'}$.

First, the four corner nodes of the oversampling domain $\Omega_{\text{os}}^{j}$ are labeled as $1', 2', 3', 4'$, and local displacement basis vectors are constructed for the x- and y-direction degrees of freedom at these four corner nodes, respectively. Taking the x-direction degree of freedom at corner node $1'$ as an example, as shown in Fig. 3, this degree of freedom is set to a unit value while the remaining seven corner node degrees of freedom are all set to zero. Linear interpolation is then performed along each edge of the outer boundary of the oversampling domain based on the displacements of the adjacent corner nodes, which serves as the displacement boundary condition for the local elasticity problem. Under the premise that the linear elasticity problem satisfies the superposition principle, the 8 local basis vectors obtained in this manner can reproduce rigid-body translational modes; when the boundary interpolation is capable of reproducing a linear displacement field, consistency of the rigid-body rotational mode can also be maintained. Therefore, the constructed oversampled numerical basis functions satisfy the necessary rigid-body-mode consistency, thereby avoiding spurious strain energy in the substructure under rigid-body motion. It should be emphasized that the aforementioned linear displacement conditions are imposed on the outer boundary of the oversampling domain $\Omega_{\text{os}}^{j}$ solely for the purpose of constructing the local auxiliary basis vectors, and do not constitute a linear assumption on the boundary displacement of the target substructure region $\Omega_{\text{sub}}^{j}$.

Under the prescribed displacement boundary conditions, finite element analysis is performed on the domain $\Omega_{\text{os}}^{j}$, which contains $(m+2l)\times(m+2l)$ fine-scale elements. This yields a displacement vector $\boldsymbol{\psi}_x^{1'} \in \mathbb{R}^{2(m+2l+1)^2}$. Similarly, setting the $y$-direction displacement at corner node $1'$ to unity and its $x$-direction displacement to zero gives $\boldsymbol{\psi}_y^{1'} \in \mathbb{R}^{2(m+2l+1)^2}$. By performing similar treatment for nodes $j' = 2, \dots, 4$ in sequence, 8 displacement basis vectors $\boldsymbol{\psi}_x^{1'}, \boldsymbol{\psi}_y^{1'}, \dots, \boldsymbol{\psi}_x^{4'}, \boldsymbol{\psi}_y^{4'}$ can be obtained.

Furthermore, taking the $x$-direction basis vectors as an example, the elements corresponding to the boundary and interior node displacements of the substructure $\Omega_{\text{sub}}^j$ are extracted from $\boldsymbol{\psi}_x^{j'}(j' = 1,2,3,4)$, yielding the corresponding target substructure displacement basis vectors $\boldsymbol{\varphi}_x^{j'} \in \mathbb{R}^{2(m+1)(m+1)}(j' = 1,2,3,4)$. By assembling $\boldsymbol{\varphi}_x^{j'}(j' = 1,2,3,4)$ column-wise into a matrix, a numerical basis function matrix corresponding to the substructure $\Omega_{\text{sub}}^j$ and its 4 corner nodes can be obtained. Similarly, the same extraction operation is performed on the y-direction basis vectors; finally, the 8-target substructure displacement basis vectors are assembled column-wise in the order shown in Eq. (10) to form the initial numerical basis function matrix:

$$\boldsymbol{\phi} = \left[\boldsymbol{\varphi}_x^{1'}\ \boldsymbol{\varphi}_y^{1'}\ \boldsymbol{\varphi}_x^{2'}\ \boldsymbol{\varphi}_y^{2'}\ \boldsymbol{\varphi}_x^{3'}\ \boldsymbol{\varphi}_y^{3'}\ \boldsymbol{\varphi}_x^{4'}\ \boldsymbol{\varphi}_y^{4'}\right] = [\boldsymbol{\phi}_1, \dots, \boldsymbol{\phi}_8], \tag{10}$$

where $\boldsymbol{\phi} \in \mathbb{R}^{2(m+1)(m+1)\times 8}$.

In fact, the 8 column vectors in $\boldsymbol{\phi}$ each represent a distinct deformation mode. Let the index set of the 8 corner node degrees of freedom of the substructure $\Omega_{\text{sub}}^j$ be denoted as $\chi_{\text{c}}$, and further define a restriction operator $\mathcal{M}$ for extracting the corresponding rows of the matrix $\boldsymbol{\phi}$ according to the corner node degree-of-freedom index $\chi_{\text{c}}$ of the target substructure. This yields the matrix $\mathbf{T} = \mathcal{M}(\boldsymbol{\phi}) \in \mathbb{R}^{8\times 8}$. In general, the matrix $\mathbf{T}$ obtained by the above procedure may not be an identity matrix (i.e., $\mathbf{T} \neq \mathbf{I}$), which would destroy the Kronecker-δ property that the numerical basis function matrix should possess. To address this, we introduce coefficients $C_{ij}, i, j = 1, \dots, 8$ to linearly combine the column vectors of the original column of the matrix $\boldsymbol{\phi}$ i.e., let:

$$\widetilde{\boldsymbol{\phi}}_i = \sum_{j=1}^{8} C_{ji}\, \boldsymbol{\phi}_j, \quad i, j = 1, \dots, 8, \tag{11}$$

Then the new interpolation shape function matrix is formed from $\widetilde{\boldsymbol{\phi}}_i$.

$$\widetilde{\boldsymbol{\phi}} = \left[\widetilde{\boldsymbol{\phi}}_1, \dots, \widetilde{\boldsymbol{\phi}}_8\right] = \boldsymbol{\phi}\mathbf{C}, \tag{12}$$

where $\mathbf{C} = \left(C_{ij}\right) \in \mathbb{R}^{8\times 8}$.

The interpolation coefficients $C_{ij}, i, j = 1, \dots, 8$ should be determined such that $\widetilde{\boldsymbol{\phi}}$ satisfies the Kronecker-δ property, i.e., it should satisfy:

$$\mathcal{M}\left(\widetilde{\boldsymbol{\phi}}\right) = \mathcal{M}(\boldsymbol{\phi}\mathbf{C}) = \mathcal{M}(\boldsymbol{\phi})\mathbf{C} = \mathbf{T}\mathbf{C} = \mathbf{I}_{8\times 8}, \tag{13}$$

Solving these yields:

$$\mathbf{C} = \mathbf{T}^{-1}. \tag{14}$$

During the above normalization process, it is necessary that the matrix $\mathbf{T}$ remains non-singular. For regular substructures with effective material distribution, non-degenerate corner node degree-of-freedom definitions, and an appropriately sized oversampling domain, $\mathbf{T}$ is numerically invertible; if ill-conditioned or near-singular situations arise in actual computations, they need to be handled through regularization, sample removal, or local exact solution methods.

Based on the steps above, an oversampled numerical basis function matrix corresponding to the substructure $\Omega_{\text{sub}}^{j}$ and its 4 corner nodes, satisfying the Kronecker-δ property, can be obtained:

$$\widetilde{\boldsymbol{\phi}} = \boldsymbol{\phi}\mathbf{C} = \boldsymbol{\phi}\mathbf{T}^{-1} = \boldsymbol{\phi}\big(\mathcal{M}(\boldsymbol{\phi})\big)^{-1}. \tag{15}$$

At this point, the oversampled numerical basis function matrix $\widetilde{\boldsymbol{\phi}}$ possesses properties similar to $\boldsymbol{N}^{j}$ in Eq. (7):

$$\boldsymbol{u}_{j} = \widetilde{\boldsymbol{\phi}}^{j}\boldsymbol{u}_{\text{c}}^{j}, \tag{16}$$

where $\widetilde{\boldsymbol{\phi}}^{j} \in \mathbb{R}^{ns\times nc}$ is the oversampled numerical basis function matrix that maps the corner node degrees of freedom of the substructure to all node degrees of freedom of that substructure, and $ns$ and $nc$ denote, respectively, the total number of degrees of freedom of all nodes in the substructure and the total number of corner node degrees of freedom. As shown in Fig. 1, for a two-dimensional problem discretized with quadrilateral elements, the nodes labeled $1, 2, 3, 4$ in red are the corner nodes, in which case $nc = 8$. The remaining nodes labeled in red are non-corner nodes. According to Eq. (16), the total strain energy of substructure j can be expressed as:

$$W^{j} = \frac{1}{2}\left(\boldsymbol{u}^{j}\right)^{\top}\mathbf{K}^{j}\boldsymbol{u}^{j} = \frac{1}{2}\left(\widetilde{\boldsymbol{\phi}}^{j}\boldsymbol{u}_{\text{c}}^{j}\right)^{\top}\mathbf{K}^{j}\widetilde{\boldsymbol{\phi}}^{j}\boldsymbol{u}_{\text{c}}^{j} = \frac{1}{2}\left(\boldsymbol{u}_{\text{c}}^{j}\right)^{\top}\mathbf{K}_{\text{r}}^{j}\boldsymbol{u}_{\text{c}}^{j}, \tag{17}$$

where $\mathbf{K}_{\text{r}}^{j}$ is the local stiffness matrix projected onto the oversampled reduced space:

$$\mathbf{K}_{\text{r}}^{j} = \left(\widetilde{\boldsymbol{\phi}}^{j}\right)^{\top}\mathbf{K}^{j}\widetilde{\boldsymbol{\phi}}^{j}. \tag{18}$$

It should be emphasized that what has been obtained so far is the oversampling reduced-order mapping and its local condensed stiffness on a single substructure. Since the numerical basis functions of each substructure are constructed independently, if they are directly stitched into a global mapping, the displacements on the shared edges of adjacent substructures are generally not automatically continuous. Therefore, this paper only gives the formal representation of the subsequent global reduced-order equations here; their rigorous construction depends on the global compatible displacement field established in Section 3.

After defining the global compatible space in Section 3, analogous to the assembly process of the global stiffness matrix $\mathbf{K}_{\text{global}=} \bigvee_{j=1}^{N} \mathbf{K}^{j}$ in finite element analysis, $\widetilde{\boldsymbol{\Phi}}^{j}$ can also be embedded into the global reduced-order mapping $\widetilde{\boldsymbol{\Phi}}_{\text{global}=} \bigvee_{j=1}^{N} \widetilde{\boldsymbol{\Phi}}^{j}$, thereby formally obtaining the condensed global stiffness matrix $\mathbf{K}_{\text{global}}^{\text{c}}$:

$$\mathbf{K}_{\text{global}}^{\text{c}} = \left(\widetilde{\boldsymbol{\Phi}}_{\text{global}}\right)^{\text{T}} \mathbf{K}_{\text{global}} \widetilde{\boldsymbol{\Phi}}_{\text{global}}. \tag{19}$$

After obtaining the globally compatible reduced-order stiffness matrix, the substructure corner node degrees of freedom can be solved through the reduced-order linear system:

$$\mathbf{K}_{\text{global}}^{\text{c}} \boldsymbol{u}_{\text{c}} = \boldsymbol{f}_{\text{eq}}. \tag{20}$$

Therefore, after the global compatibility mechanism is established, the solution scale of Eq. (20) can be reduced from the fine-scale finite element nodal degrees of freedom to the substructure corner node degrees of freedom. After obtaining the corner node displacements, the fine-scale displacements within each substructure can be recovered using Eq. (16). This recovery process is essentially a matrix multiplication, and thus typically incurs very low computational overhead relative to solving the global linear system.

Let Ω denote the region occupied by the global structure requiring finite element analysis. For substructures adjacent to the boundary $\partial\Omega$ of the global structure, the outward-extended layer $\Omega_{\text{os}}^{j} \backslash \Omega_{\text{sub}}^{j}$ of the oversampling domain relative to the substructure region will partially extend beyond Ω. In this case, "virtual" material parameters $\boldsymbol{E}_{\text{vir}}^{j}$ need to be specified for the extended portion in order to compute the oversampled numerical basis function $\widetilde{\boldsymbol{\Phi}}$. In this paper, a classified extrapolation strategy is adopted, as shown in Fig. 4.

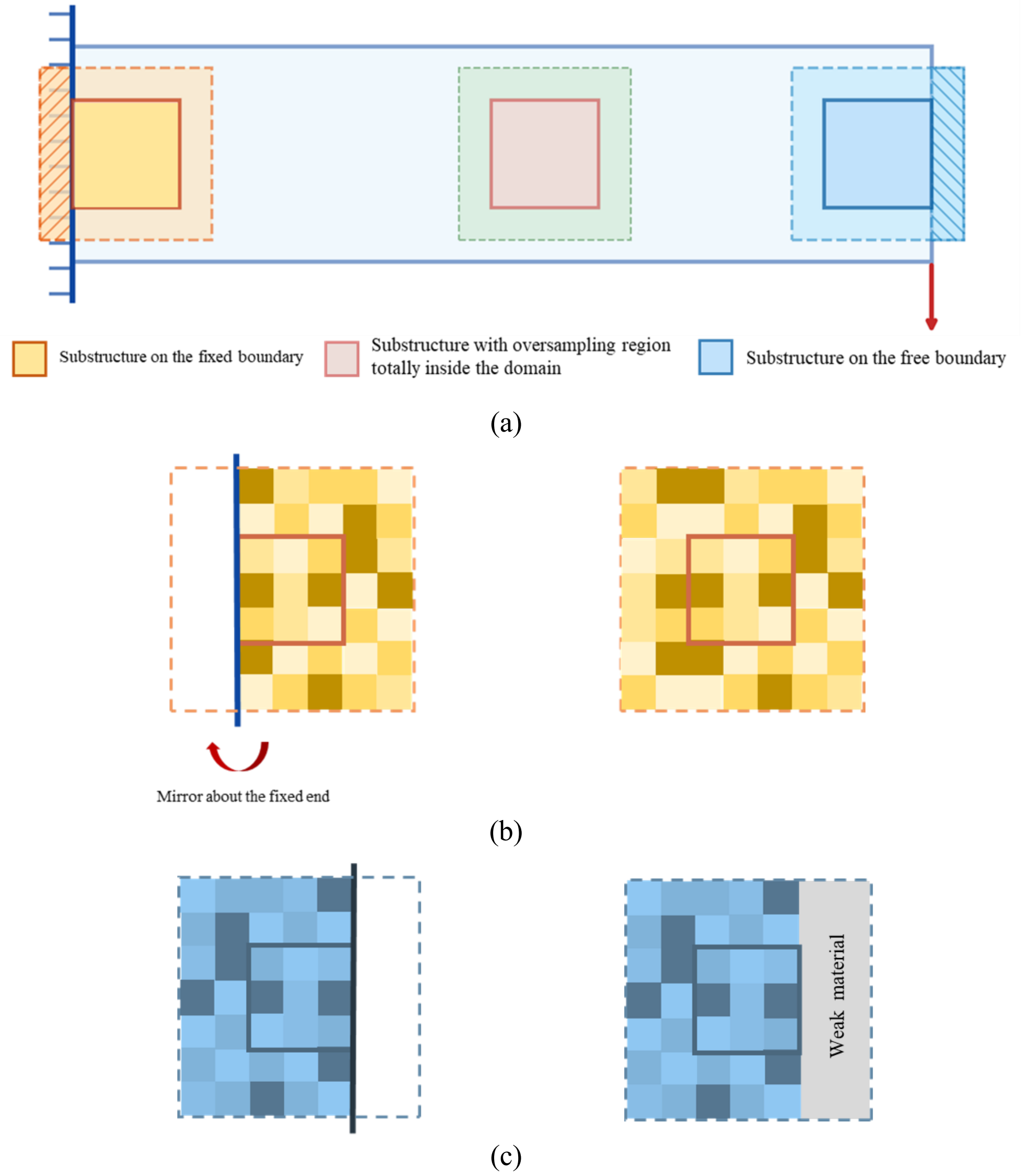

Fig. 4. Treatment of oversampling domains that extend beyond the physical structural boundary. (a) Classification of boundary-adjacent substructures. (b) Material extension across a fixed boundary. (c) Material extension across a free boundary.

For extensions beyond the free-boundary direction (edges without any Dirichlet constraint), the material parameters of all extended elements are set identical to those of the low-stiffness material. For extensions beyond the fixed-boundary direction (or boundaries where Dirichlet constraints are imposed on the global structure), a mirror-symmetric extrapolation with respect to the structural boundary is employed. That is, the virtual material parameter at the extended point x is taken as the true material parameter at its mirror-image point, denoted as $\boldsymbol{E}_{\mathrm{vir}}(\boldsymbol{x}) = \boldsymbol{E}(\boldsymbol{\mathcal{R}}_{\Gamma}(\boldsymbol{x}))$, where $\boldsymbol{\mathcal{R}}_{\Gamma}$ is the mirror mapping of the outward-extended layer with respect to the fixed end face.

This mirroring operation is used only to construct the virtual material neighborhood outside the fixed boundary; the true Dirichlet condition is still imposed in the global analysis, and should not be interpreted as an additional even-symmetry assumption imposed on the true displacement field.

The aforementioned virtual material extrapolation is used only for substructures near the global boundary when constructing the oversampled numerical basis functions, in order to supplement the extended region; it does not alter the true material distribution within the real structural domain Ω, nor does it replace the true boundary conditions imposed in the global finite element analysis.

# 3 Overlapping Finite Element Displacement Compatibility Method Based on Partition of Unity

In Section 2, an oversampled numerical basis function $\widetilde{\boldsymbol{\Phi}}^j$ was constructed for each substructure, as shown in Eq. (16). This mapping is used to recover the fine-scale displacements within the substructure from its corner node degrees of freedom. Since $\widetilde{\boldsymbol{\Phi}}^j$ satisfies the Kronecker-δ property, after global assembly at the shared corner node degrees of freedom, adjacent substructures can maintain consistent displacements at the shared corner nodes. However, at non-corner nodes on the shared edge, the displacements on the two sides are determined by their respective oversampling domains and local material distributions, and are generally not the same function. As a result, the following situation may occur:

$$\boldsymbol{u}_{j_1}(\boldsymbol{x}) \neq \boldsymbol{u}_{j_2}(\boldsymbol{x}), x \in \left(\partial\Omega_{\mathrm{sub}}^{j_1} \cap \partial\Omega_{\mathrm{sub}}^{j_2}\right) \backslash \Omega_{\mathrm{c}}^{j_1 j_2} \tag{21}$$

where $\boldsymbol{u}_{j_1}$ and $\boldsymbol{u}_{j_2}$ denote the displacement fields recovered by the local numerical basis functions of substructures $j_1$ and $j_2$, respectively; $\partial\Omega_{\mathrm{sub}}^{j_1}$ and $\partial\Omega_{\mathrm{sub}}^{j_2}$ are the boundaries of the two substructures; and $\Omega_{\mathrm{c}}^{j_1 j_2}$ denotes the region corresponding to the shared corner node degrees of freedom of the two substructures. Eq. (21) indicates that if the independently constructed displacements of each substructure are directly assembled, discontinuities in displacement may arise on the shared edges, such that the assumed displacement after assembly no longer belongs to the admissible displacement space required by the global weak form, thereby destroying the consistency of the global finite element discretization. The direct consequence of this is that unreliable strain jumps and unreliable strain energy calculations may occur at the interface, which may further lead to bias in topology optimization sensitivity.

Therefore, although the oversampling reduced-order mapping obtained in Section 2 is capable of describing the fine-scale displacement within a single substructure, it does not yet resolve the global compatibility problem among adjacent substructures. The focus of this section is no longer to reconstruct local numerical basis functions, but rather to establish, under the premise of maintaining a low-dimensional representation of the corner node degrees of freedom, a compatible displacement space that can be used for global assembly. To this end, this paper introduces an overlapping finite element concept based on partition of unity: first, an overlapping substructure network with shared coverage regions is constructed; then, partition-of-unity weighting functions are used to perform weighted combination of the local displacements within the overlapping regions, thereby obtaining a globally compatible displacement field that can be used for reduced-order analysis [40-48].

### 3.1 Overlapping Coverage of Substructures and Fine-Scale Mesh Construction

In order to use partition of unity in Section 3.2 to reconcile the displacements of adjacent substructures, this section first constructs an overlapping substructure coverage. The key here is not simply to introduce an overlapping region, but to ensure that the overall fine-scale mesh, after local overlap, still maintains a regular Cartesian form. To achieve this, rather than directly overlapping uniform substructure meshes, this paper first designs non-uniform boundary elements within each substructure, and then allows adjacent substructures to overlap near their shared edges according to a globally uniform element size. This provides a shared coverage region for overlapping finite element reconciliation of local displacements, while also avoiding the introduction of additional non-uniform global mesh caused by the substructure partitioning itself.

If each substructure still uses a regular, uniform $h \times h$ mesh, and adjacent substructures are made to overlap by a non-integer element width — for example, an overlap of $h/2$ — the overlapping region will be cut into non-uniform sub-elements such as $(h/2) \times (h/2)$ or $(h/2) \times h$, as shown in Fig. 5. Such small elements are not driven by displacement, stress, or material interface features, but are artificially generated by the overlapping partitioning strategy; they increase the number of global mesh element types, complicate the assembly rules, and weaken the ability of a regular mesh to directly represent material distribution. It should be noted that the overlap width also cannot simply be taken as one complete fine element, otherwise the normalization denominator of the subsequent partition-of-unity weighting function may become zero in local regions—related discussion is provided in [40-48]. Therefore, this paper does not directly overlap the uniform substructure meshes.

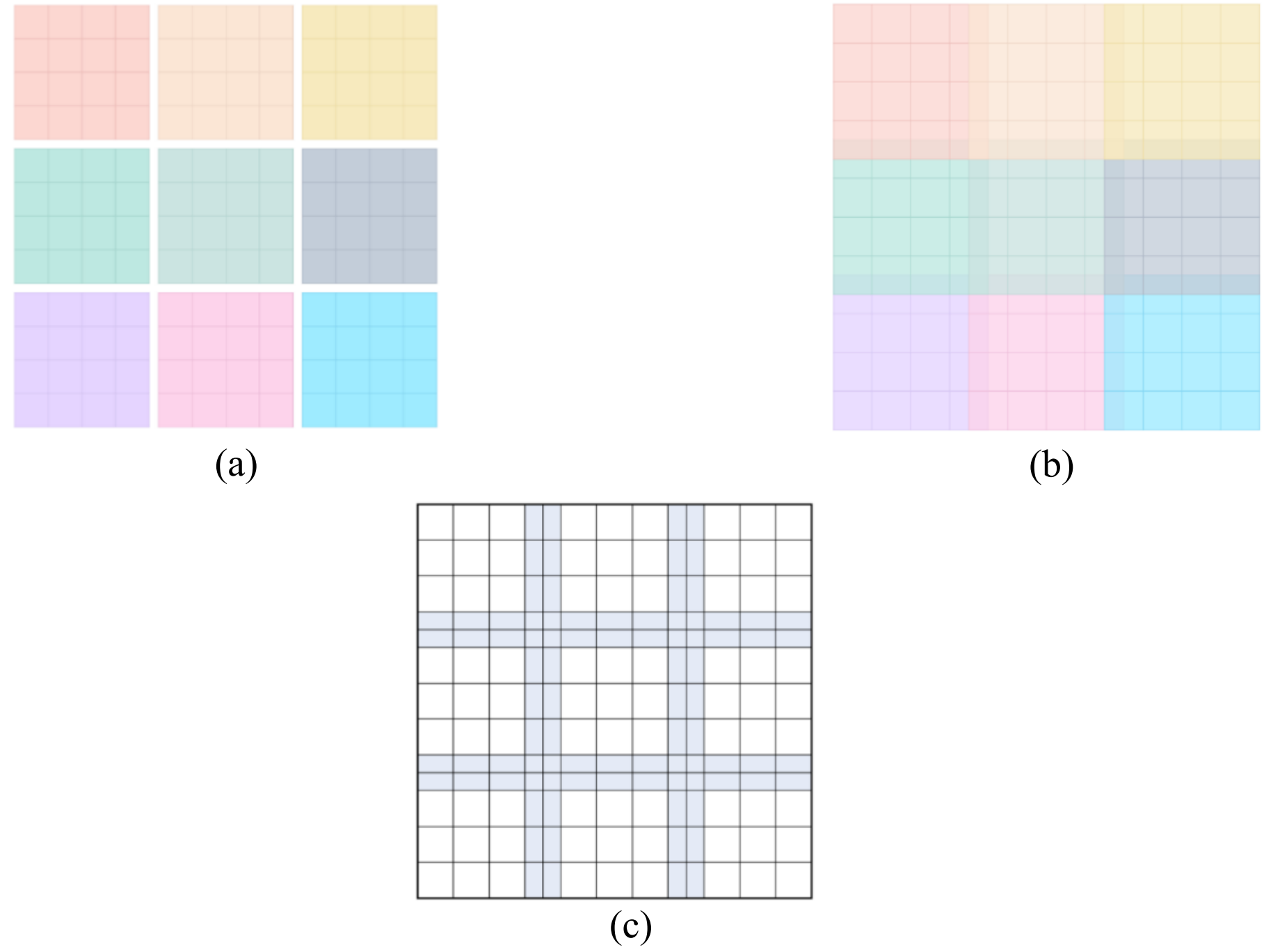

Fig.5. Nonuniform global grid produced by directly overlapping uniformly discretized substructures. (a) Nine substructures with uniform local fine-scale meshes. (b) Superposition of the nine substructures. (c) Resulting nonuniform integration grid.

To avoid the aforementioned artificial refinement, this paper first adopts a specific non-uniform mesh within each individual substructure, as shown in Fig. 6. Specifically, on the side shared with a neighbor, the outermost layer of boundary elements is expanded from $h \times h$ to $2h \times h$, with the $2h$ direction oriented perpendicular to the overlap band; if two adjacent edges at a corner node both need to overlap with neighbors, then the corresponding corner node element is simultaneously expanded in both directions, with dimensions $2h \times 2h$. In this way, after adjacent substructures overlap by a width of $h$, the overlapping region can be reassembled into a globally uniform $h \times h$ fine mesh.

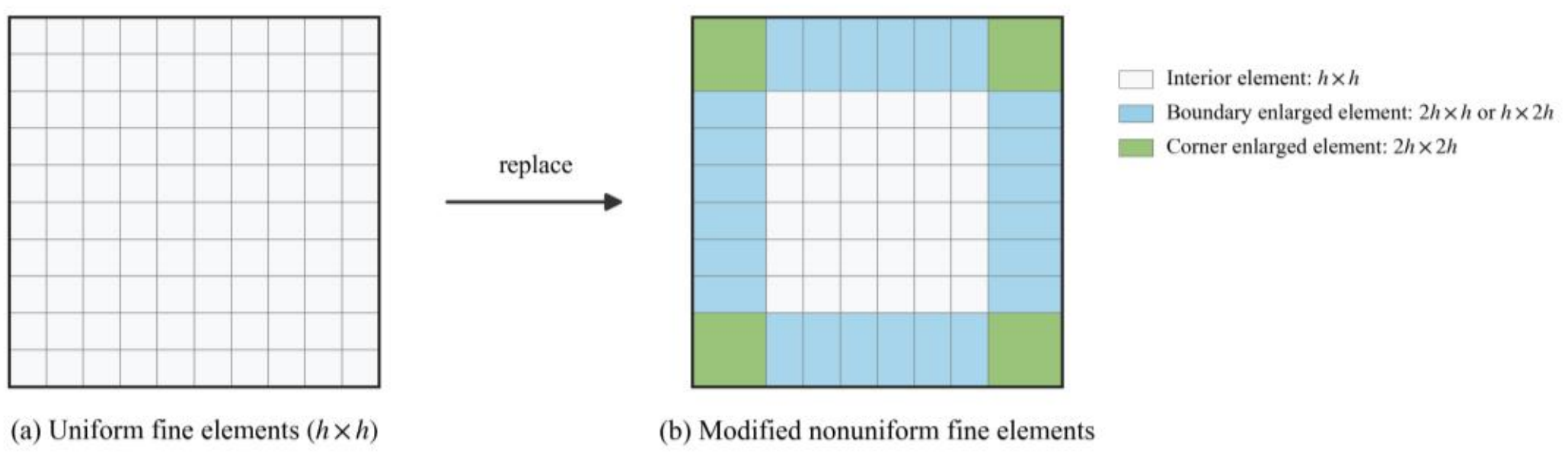

Fig.6. Replacing the uniform fine elements composing a substructure with specific non-uniform fine

elements.

Accordingly, the local mesh of an individual substructure can be classified into three types of elements: interior $h \times h$ elements, boundary $2h \times h$ expanded elements, and corner node $2h \times 2h$ expanded elements, as shown in Table 1.

Table 1. Local element categories in an overlapping substructure.

| Element category | Dimensions | Location within the substructure |
|---|---|---|
| Interior fine-scale element | $h \times h$ | Substructure interior $(m-4) \times (m-4)$ (white cells in Fig. 6) |
| Edge-extended element | $2h \times h$ or $h \times 2h$ | Overlap edge, excluding the two corner elements (blue cells in Fig. 6) |
| Corner-extended element | $2h \times 2h$ | Intersection of two orthogonal overlap bands (green cells in Fig. 6) |

Since expansion is only required on the side where an adjacent substructure exists, the type of a given substructure is determined by its position within the $n_{\mathrm{s}} = n_{\mathrm{slx}} \times n_{\mathrm{sly}}$ rectangular arrangement, resulting in a total of 9 categories, as shown in Fig. 7 and Table 2.

Table 2. Nine substructure types in a rectangular arrangement.

| Type No. | Type | Overlapping directions with adjacent substructures | Count |
|---|---|---|---|
| Type 1 | Up Left Corner | Bottom, Right | 1 |
| Type 2 | Up | Bottom, Left, Right | $n_{\mathrm{slx}} - 2$ |
| Type 3 | Up Right Corner | Bottom, Left | 1 |
| Type 4 | Left | Top, Bottom, Right | $n_{\mathrm{sly}} - 2$ |
| Type 5 | Internal | Top, Bottom, Left, Right | $(n_{\mathrm{slx}} - 2) \times (n_{sly} - 2)$ |
| Type 6 | Right | Top, Bottom, Left | $n_{\mathrm{sly}} - 2$ |
| Type 7 | Bottom Left Corner | Top, Right | 1 |
| Type 8 | Bottom | Top, Left, Right | $n_{\mathrm{slx}} - 2$ |
| Type 9 | Bottom Right Corner | Top, Left | 1 |

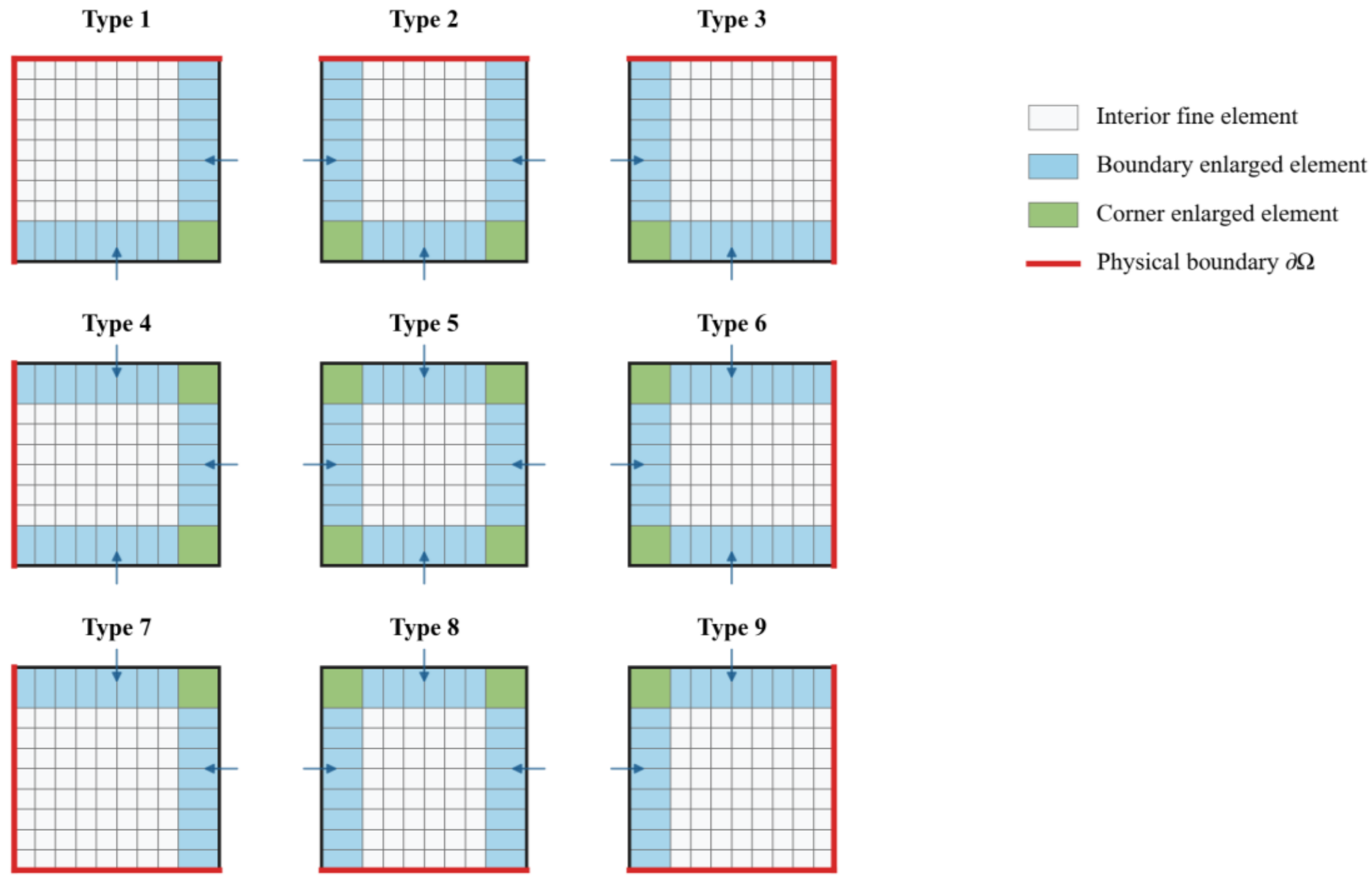

Fig. 7. Nine different types of substructures.

Specifically, the "internal" type has neighbors on all four sides, and expansion is required on all four sides; types 2, 4, 6, and 8 have one edge located on the global boundary, and expansion is not performed on that side; types 1, 3, 7, and 9 have two adjacent edges located on the global boundary, and expansion is not performed on either of those two sides. Following this rule, uniform fine-mesh discretization can be obtained after overlapping, as shown in Fig. 7 and Fig.8.

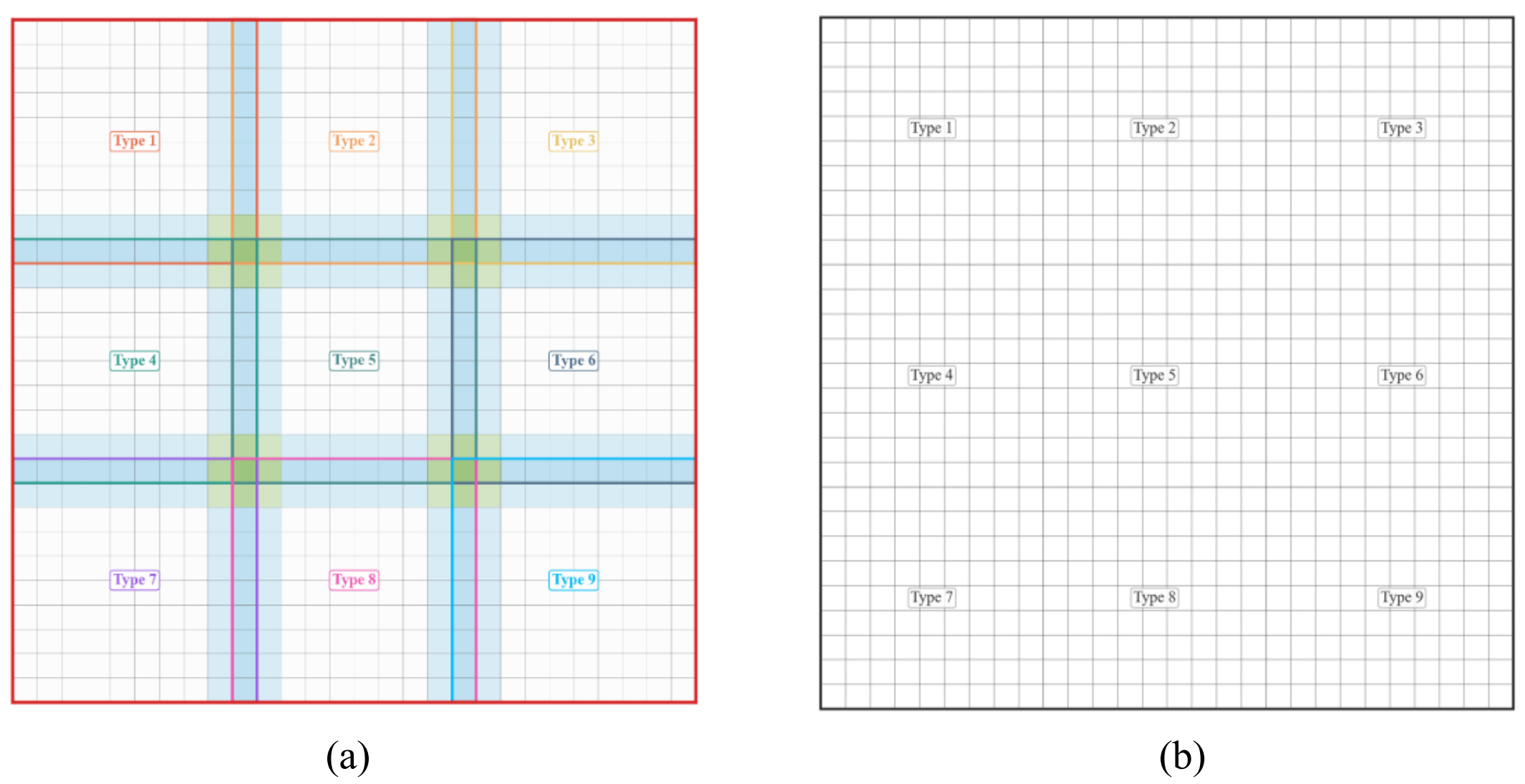

(a) (b)

Fig.8. A 3×3 arrangement containing one representative of each substructure type; their superposition produces a uniform Cartesian fine-scale grid . (a) Overlap of the nine substructure

types. (b) Resulting uniform Cartesian fine grid.

## 3.2 Displacement Compatibility Based on Overlapping Finite Elements

Section 3.1 has already given the geometric construction of the overlapping substructure mesh. It should be further clarified that, after overlapping, a global $h \times h$ fine-scale element is not an independent original element within a single substructure, but rather an accumulated region formed by the mutual overlap of local elements from adjacent substructures. Therefore, at the same physical location, there simultaneously exist interpolation expressions from different local elements; if calculated separately in the traditional finite element manner and the local stiffness matrices are directly assembled, the displacement degrees of freedom in the overlapping region will lack a unique compatibility relationship, causing the global analysis to fail. The overlapping finite element work of Bathe et al. [40-48] provides a feasible approach to this compatibility problem by using partition-of-unity weighting functions to perform weighted combination of the local interpolations.

The goal of this section is: while preserving the local element interpolation expressions, to construct a set of non-negative, locally supported weighting functions that sum to one, so that multiple local displacement fields within the overlapping region can be combined into a unique global displacement field. Specifically, let $\Omega_l$ denote the coverage region of the $l$-th local element, and $N_{\mathrm{e}}$ denote the total number of local elements. The weighting function $\omega_l(\boldsymbol{x})$ needs to satisfy the following conditions:

1) It is continuous and non-negative over the analysis domain;
2) It vanishes outside $\Omega_l$, i.e., each local interpolation only takes effect within its own coverage region;
3) It satisfies $\sum_{l=1}^{N_{\mathrm{e}}} \omega_l(\boldsymbol{x}) = 1$, over the entire domain $\Omega$.

The above conditions guarantee that: if each local interpolation is compatible within its own coverage region and can reproduce a linear displacement field, then the global interpolation obtained by partition-of-unity weighting also preserves compatibility and linear completeness, and can therefore pass the patch test.

Let

$$P_l(\boldsymbol{x}) = \sum_{r} h_r^l(\boldsymbol{x})\ p_r^l, \tag{22}$$

be a function defined on $\Omega_l$, where $h_r^l(\boldsymbol{x})$ is the standard Q4 finite element shape function, and $p_r^l$

is the marker value assigned to node r of the $l$-th local element. The role of $p_r^l$ is to identify whether that node is located on the shared coverage boundary between the current element and other elements: if the node lies at the intersection of the coverage regions of $\Omega_l$ and other local elements, then $p_r^l$ is taken as 0; otherwise $p_r^l$ is taken as 1 (including nodes located on the boundary $\partial\Omega$ of the global structure), as shown in Fig. 9. The $P_l(\boldsymbol{x})$ constructed in this way is zero on the coverage boundary, and can be continuously extended to the entire analysis domain.

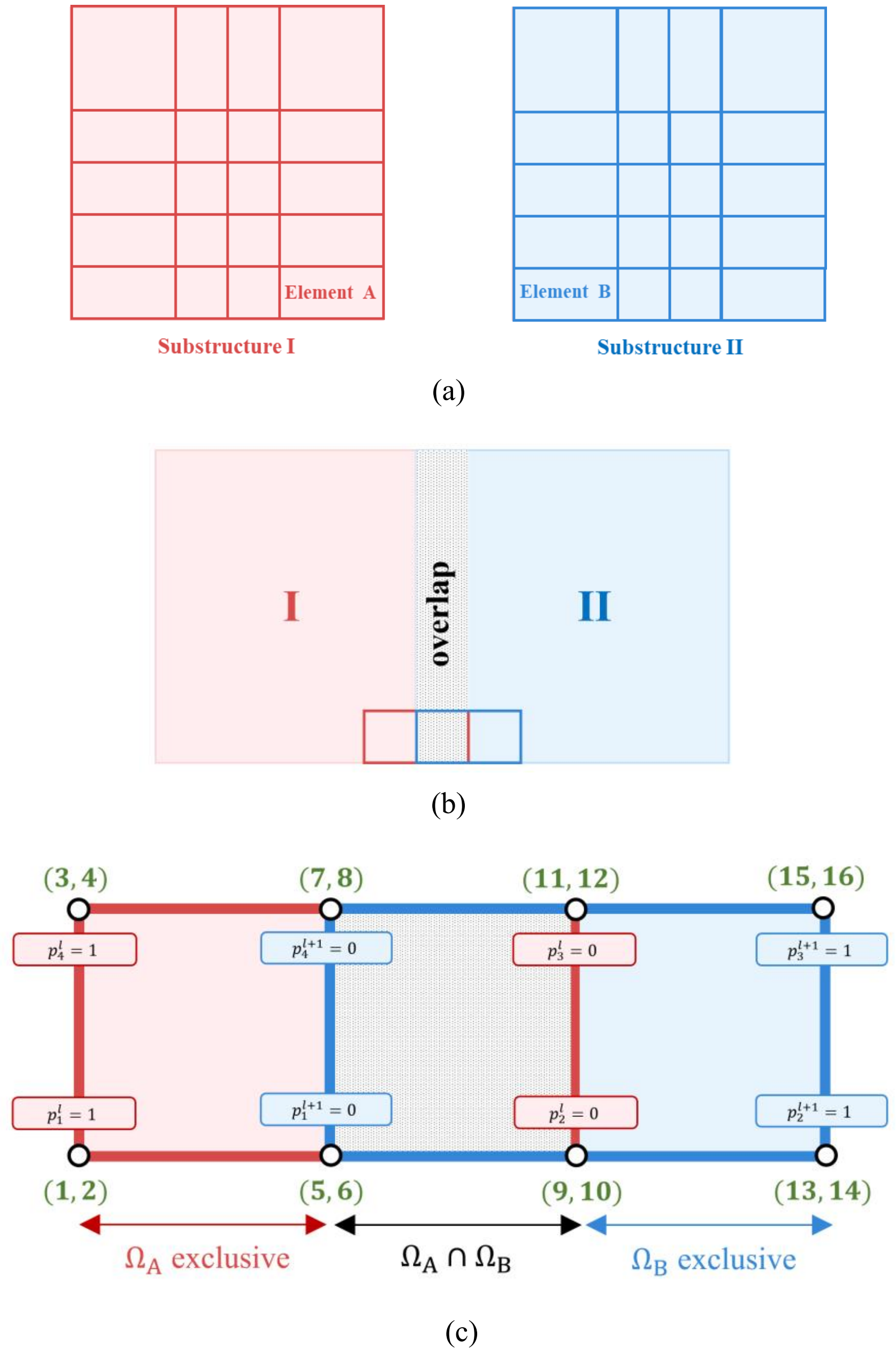

Fig.9. Schematic diagram of two $h \times 2h$-sized elements overlapping in the x-direction, with the gray region indicating the overlapping area. (a) Element A and Element B are the elements on the boundaries of Substructure I and Substructure II. (b) The positions of Element A and Element B

when Substructure I and Substructure II overlap, with the two elements covering each other. (c) The $p_r^l$ values corresponding to the nodes in the figure have been marked. The green numbers denote the global indices corresponding to each nodal degree of freedom ($x$-direction first, $y$-direction second).

To satisfy the partition-of-unity property, the weighting function is defined as:

$$w_l(\boldsymbol{x}) = \frac{P_l(\boldsymbol{x})}{\sum_l P_l(\boldsymbol{x})}. \tag{23}$$

Eq. (23) normalizes $P_l(\boldsymbol{x})$ over all local elements covering the current point, thereby obtaining a weighting function that satisfies the partition-of-unity property. For a reasonable overlapping mesh, any point within the analysis domain is covered by at least one local element satisfying $P_l(\boldsymbol{x}) > 0$, so the denominator $\sum_l P_l(\boldsymbol{x})$ is nonzero, and Eq. (23) is well-defined. [45] presents a case of an ineffective mesh, in which a local region may encounter the problem of a zero denominator.

Suppose the local interpolation of displacement has already been established on each mesh:

$$\boldsymbol{u}_l(\boldsymbol{x}) = \mathbf{H}_l(\boldsymbol{x})\boldsymbol{q}_l, \tag{24}$$

where $\boldsymbol{u}_l(\boldsymbol{x})$ denotes the interpolated displacement vector on the $l$-th local element, $\mathbf{H}_l(\boldsymbol{x})$ denotes the shape function matrix of that element, and $\boldsymbol{q}_l$ denotes the corresponding degree-of-freedom vector. Since partition of unity is only responsible for global compatibility, the local interpolation itself may adopt any compatible interpolation technique; this paper adopts the traditional isoparametric interpolation.

Thus, the global displacement can be written as the partition-of-unity weighted sum of local displacements:

$$\boldsymbol{u}_{\text{global}}(\boldsymbol{x}) = \sum_{l=1}^{\boldsymbol{N}_\text{e}} w_l(\boldsymbol{x})\,\boldsymbol{u}_l(\boldsymbol{x}) = \sum_{l=1}^{\boldsymbol{N}_\text{e}} w_l(\boldsymbol{x})\,\mathbf{H}_l(\boldsymbol{x})\boldsymbol{q}_l, \tag{25}$$

The strain interpolation can be directly obtained from the displacement:

$$\boldsymbol{\varepsilon}(\boldsymbol{x}) = [\mathbf{B}_1 \quad \cdots \quad \mathbf{B}_l]\begin{bmatrix}\boldsymbol{q}_1 \\ \vdots \\ \boldsymbol{q}_l\end{bmatrix} = \sum_{l=1}^{\boldsymbol{N}_\text{e}} \mathbf{B}_l(\boldsymbol{x})\,\boldsymbol{q}_l, \tag{26}$$

where

$$\mathbf{B}_l = \begin{bmatrix} \frac{\partial}{\partial x} & 0 \\ 0 & \frac{\partial}{\partial y} \\ \frac{\partial}{\partial y} & \frac{\partial}{\partial x} \end{bmatrix} w_l(\boldsymbol{x})\mathbf{H}_l(\boldsymbol{x}), \tag{27}$$

The stiffness matrix is written as

$$\mathbf{K}_{\text{global}} = \int_{\Omega} [\mathbf{B}_1 \quad \cdots \quad \mathbf{B}_l]^{\top} \mathbf{C}\, [\mathbf{B}_1 \quad \cdots \quad \mathbf{B}_l] \mathrm{d}\Omega, \tag{28}$$

where $\mathbf{C}$ is the elastic stress-strain matrix. Due to the local support property of the weight functions, coupling between different local elements only occurs within their overlapping regions. Therefore, for two overlapping fine elements $\alpha$ and $\beta$ , the corresponding coupling stiffness term can be defined:

$$\mathbf{K}_{\text{global}}(\alpha, \beta) = \int_{\Omega} \mathbf{B}_{\alpha}{}^{\top} \boldsymbol{C}\, \mathbf{B}_{\beta} \mathrm{d}\Omega = \int_{\Omega_{\alpha} \cap \Omega_{\beta}} \mathbf{B}_{\alpha}{}^{\top} \boldsymbol{C}\, \mathbf{B}_{\beta} \mathrm{d}\Omega. \tag{29}$$

The last equality in Eq. (29) holds because $\mathbf{B}_{\alpha}$ and $\mathbf{B}_{\beta}$ vanish outside $\Omega_{\alpha}$ and $\Omega_{\beta}$, respectively. This indicates that the key to overlapping finite element assembly is not simply summing the local stiffness matrices, but rather integrating the coupling terms over the overlapping region for piecewise continuity. Taking Fig. 9 as an example, after the red and blue local elements overlap, their union can be divided into three integration sub-regions: the region exclusively occupied by the red element, the mutual overlapping region of the two, and the region exclusively occupied by the blue element. The middle overlapping region is not an independent physical element in the original discretization, but rather an integration sub-region introduced for numerical integration and stiffness assembly. In implementation, these integration sub-regions can be uniformly treated as global $h \times h$ fine-scale elements, and the stiffness contributions are assembled according to the local element degrees of freedom associated with them. Therefore, at the level of stiffness matrix assembly, the overall structure $\Omega$ can be equivalently regarded as being covered by $(n_{\text{slx}}(m-1)+1) \times (n_{\text{sly}}(m-1)+1)$ $h \times h$-sized fine elements.

Unlike general overlapping finite element methods [40-48], the integration and assembly in this work are ultimately performed on a uniform Cartesian fine-scale mesh; therefore, all $h \times h$ elements are congruent in the reference coordinates $(\xi, \eta) \in [-1,1]^2$ , with $\mathbf{J} = (h/2)\mathbf{I}$ and $\det(\mathbf{J}) = (h^2/4)$.

This means that all $h \times h$ elements possess three properties: 1) Constancy, i.e., $\mathbf{J}_e$ does not depend on $(x, y)$; 2) Uniformity, i.e., $\mathbf{J}_e$ does not depend on the element index e; 3) Isotropy and orthogonality, i.e., $\mathbf{J}_e$ is a scalar multiple of the identity matrix, so that the x and y directions have the same stretching ratio. It follows that, as long as two integration sub-regions have the same weighting function combination form in reference coordinates, their integrand forms are identical, and the corresponding unit Young's modulus stiffness matrices are also identical. Therefore, in the method of this paper, geometric differences are all eliminated by the unified Cartesian fine mesh, and the remaining differences arise mainly from the weighting function types determined by the local element overlap relationships.

Specifically for two-dimensional problems, all fine elements can be classified into at most 16 types, based on the number and relative positions of local elements simultaneously covering each fine element, as detailed in Table 3 and Fig. 10. For each type, the corresponding integrand has a fixed polynomial form in reference coordinates and is independent of the material distribution. This yields two implementation advantages:

**1) Numerical integration can be precomputed offline and stored for reuse.** For each type $t \in \{1,2, \ldots, 16\}$, the integrand is a low-order polynomial in $(\xi, \eta)$. Therefore, 3-point Gauss quadrature in each reference coordinate direction (i.e., a $3 \times 3$ Gauss rule on $[-1,1]^2$) is adopted. Since this process does not depend on the element Young's modulus of the $e$-th element $E_e$, for analysis or topology optimization processes that require repeatedly reassembling the stiffness matrix, the 16 unit-Young's modulus stiffness matrices $\mathbf{K}_0^t$ $(t = 1,2, \ldots, 16)$ can be computed once before the main loop begins and stored in memory for reuse; subsequently, each assembly only requires a single scalar scaling of the current element:

$$\mathbf{K}_e = E_e \mathbf{K}_0^t \tag{30}$$

**2) The assembly process is simplified.** Due to the regularity of the Cartesian mesh, the global DOF indices of each fine element can be obtained directly from the mesh index. Combined with the 16 types of precomputed unit stiffness matrices $\mathbf{K}_0^t$, the assembly of the entire global stiffness matrix $\mathbf{K}_{\text{global}}$ can be completed through the workflow "grouping by type — vectorized scaling by Young's modulus — assembly via predetermined indices," and can be obtained by the single sparse matrix assembly. The runtime stage mainly consists of two steps—Young's-modulus scaling and sparse assembly pass-through sparse matrix assembly. For a fixed local overlap pattern and the finite set of 16 element types, the assembly cost scales linearly with the number of fine-scale integration elements, and elements of different types are mutually independent, allowing natural parallelization.

In summary, although the overlapping finite element discretization framework of this paper theoretically introduces piecewise continuous integrands and inter-element coupling terms, since the final integration domain is organized as a unified Cartesian fine mesh, all complex integrals related to the weighting functions can be precomputed offline and reused. The runtime stage mainly consists of two steps—material modulus scaling and sparse assembly—so obtaining the global compatible stiffness matrix $\mathbf{K}_{\text{global}}$ does not change the asymptotic complexity of online assembly. After obtaining $\mathbf{K}_{\text{global}}$, it is further necessary to embed the substructure oversampling reduced-order mapping obtained in Section 2 into this compatible displacement space. The specific coupling relationship between the two and the construction of the condensed stiffness will be given in Section 3.3.

Table 3. Fine-scale element types and the number of covering substructures, $N_{\text{OV}}$ (see Fig. 10).

| Type $t$ | $N_{\text{OV}}$ | Dimension of the local stiffness block | Location |
| --- | --- | --- | --- |
| 1 | 1 | $8 \times 8$ | Fine element of size $h \times h$ inside the substructure |
| 2–5 | 1 | $8 \times 8$ | Left/right/upper/lower half of the edge overlap region in the $x/y$ direction (4 sub-types) |
| 8-11 | 1 | $8 \times 8$ | Quarter of the substructure adjacent to the corner overlap region (top-left / top-right / bottom-left / bottom-right) |
| 6,7,12,13,14,15 | 2 | $16 \times 16$ | Edge overlap zone covered by two substructures |
| 16 | 4 | $32 \times 32$ | Corner overlap zone covered by four substructures |

## 3.3 Coupling of the Oversampling Reduced-Order Mapping with Overlapping Finite Elements

Section 2 constructed the oversampled numerical basis functions on a single substructure, while Sections 3.1 and 3.2 gave, respectively, the substructure coverage scheme that can overlap into a uniform Cartesian fine mesh, and the global displacement compatibility and stiffness assembly method based on partition-of-unity weighting functions. At this point, the local reduced-order representation and the global compatibility space have each been established separately, but the two still need to be explicitly connected within the same analysis framework: if only the oversampled shape functions are used without introducing overlapping compatibility, the displacements at shared edges of adjacent substructures are generally not continuous; if only overlapping finite elements are used without introducing oversampling reduction, the global unknowns remain at the fine-scale degree-of-freedom level, making it difficult to obtain the reduction effect expected from the substructure method. Therefore, this section discusses how to embed the oversampling reduced-order mapping into the compatible displacement space defined by the overlapping finite elements, thereby obtaining a reduced-order computational model that no longer needs to directly prescribe linear or higher-order interpolation forms for the target substructure boundary displacements, while still maintaining global displacement continuity.

From the perspective of the problem structure, the oversampling method and the overlapping finite elements act on different scales. For any substructure$\Omega_{\mathrm{sub}}^{j}$, the oversampling procedure first solves a local boundary value problem, based on the local material distribution, over the extended region $\Omega_{\mathrm{os}}^{j}$ covering that substructure, and obtains, according to Eq. (16), a mapping from the corner node degrees of freedom $\boldsymbol{u}_{\mathrm{c}}^{j}$ to the substructure's fine-scale degrees of freedom $\boldsymbol{u}_{j}$. The oversampled numerical basis function $\widetilde{\boldsymbol{\Phi}}^{j}$ therein does not require the target substructure's boundary displacement to be obtained through linear interpolation of the corner node displacements, and can therefore more fully reflect local heterogeneous effects such as perturbations near voids or abrupt changes at material interfaces. At the same time, once the substructure type, oversampling range, and boundary treatment method are given, $\widetilde{\boldsymbol{\Phi}}^{j}$is mainly determined by the material distribution within the oversampling region, and its construction does not take any specific global loading condition or global displacement solution as input; this locality provides a basis for the subsequent introduction of problem-independent machine learning prediction. However, Eq. (16) still only gives a local representation within a single substructure. When multiple substructures are juxtaposed, although adjacent substructures can maintain consistent displacement at shared corner nodes through global assembly, at non-corner positions on the shared edges, the displacements on the two sides are each determined by their own oversampling region and local material distribution, and generally do not correspond to the same function. In other words, what the oversampled shape functions solve is the problem of how local complex deformations can be represented by low-dimensional corner node degrees of freedom, but it has not yet been automatically guaranteed that these local displacement fields can be assembled into an admissible displacement space satisfying the weak-form compatibility requirements of the global problem.

Overlapping finite elements precisely fill in this missing link. Section 3.2 has already explained that, once the overlapping coverage is formed, the same physical location may simultaneously have interpolated expressions coming from different local elements; the role of the partition-of-unity weighting functions $w_l(\boldsymbol{x})$ is exactly to weight-combine these local displacements into a single, unique global displacement field. When this idea is combined with the oversampling reduction, the key does not lie in redefining the local shape functions, but rather in changing the source of the local displacement degrees of freedom: in the traditional overlapping finite element method, the local element degrees of freedom $\boldsymbol{q}_l$ directly participate in the global solution; whereas in the framework of this paper, $\boldsymbol{q}_l$ is no longer retained as an independent unknown, but is instead first recovered as fine-scale displacement through the oversampling mapping from the corner node degrees of freedom of the substructure to which it belongs, and then extracted as the local element's degrees of freedom through a restriction operator. Let $R_l(\cdot)$ denote the restriction operator that extracts the $l$-th local element's degrees of freedom from the fine-scale displacement of the substructure to which it belongs; then the global displacement can be written as:

$$\boldsymbol{u}_{\text{global}}(\boldsymbol{x}) = \sum_{l} \omega_l(\boldsymbol{x})\mathbf{H}_l(x)\boldsymbol{q}_l, \qquad \boldsymbol{q}_l = R_l\left(\widetilde{\boldsymbol{\Phi}}^j \boldsymbol{u}_{\text{c}}^j\right) \tag{31}$$

Eq. (31) gives the division of labor in the coupling chain described above. First, the oversampling mapping is completed by $\widetilde{\boldsymbol{\Phi}}^j$, which expands the eight corner node degrees of freedom of each substructure into the fine-scale displacements within the region covered by that substructure; subsequently, the restriction operator extracts the displacement degrees of freedom of the corresponding local element, while the partition-of-unity weighting function $\omega_l$ weights multiple local displacements within the overlapping region, eliminating displacement jumps at shared edges. The resulting $\boldsymbol{u}_{\text{global}}(\boldsymbol{x})$ thus retains the expressive power of the oversampled shape functions for complex local displacement patterns, while also satisfying the global continuity and linear completeness guaranteed by partition-of-unity interpolation. The corner node degrees of freedom serve throughout this coupling process as the shared reduced-order control variables: they are both the input to the oversampled shape functions and the node quantities that adjacent substructures naturally share geometrically. Precisely for this reason, this paper is able to simultaneously achieve local accuracy improvement and global displacement compatibility without increasing the boundary control degrees of freedom.

At the discrete stiffness level, the above coupling corresponds to a standard Galerkin projection. Let the global fine-scale stiffness matrix obtained from the overlapping finite element assembly in Section 3.2 be $\mathbf{K}_{\mathrm{global}}^{\mathrm{OV}}$. Assembling each substructure's oversampling mapping according to the global numbering of the local fine-scale degrees of freedom and corner node degrees of freedom yields the global reduced-order mapping:

$$\widetilde{\mathbf{\Phi}}_{\mathrm{global}}= \bigvee_{j=1}^{N} \widetilde{\mathbf{\Phi}}^{j} . \tag{32}$$

Then the condensed reduced-order stiffness matrix is:

$$\mathbf{K}_{\mathrm{global}}^{\mathrm{c,OV}} = \left(\widetilde{\mathbf{\Phi}}_{\mathrm{global}}\right)^{\mathrm{T}} \mathbf{K}_{\mathrm{global}}^{\mathrm{OV}} \widetilde{\mathbf{\Phi}}_{\mathrm{global}}. \tag{33}$$

Eq. (33) is formally consistent with the global condensation projection given in Eq. (19) of Section 2, but its substantive content has changed due to the introduction of overlapping finite elements. Here, $\mathbf{K}_{\mathrm{global}}^{\mathrm{OV}}$ already includes the piecewise-continuous coupling stiffness contributions over the overlapping regions, and is no longer a simple assembly of each substructure's independent fine mesh stiffness; correspondingly, $\mathbf{K}_{\mathrm{global}}^{\mathrm{c,OV}}$ is also not a direct summation of each substructure's local condensed stiffness, but rather the projection result of the oversampling reduced-order basis in the global compatible displacement space. It is precisely this difference that allows the reduced-order system, after solving for the corner node degrees of freedom, to recover a globally continuous displacement field usable for strain, strain-energy, and topology optimization sensitivity computations. If, without introducing overlapping compatibility, the independently recovered displacements of each substructure were simply patched together, non-physical jumps would often appear in the interfacial strain field, and the subsequent sensitivity analysis based on strain energy density would also lose validity as a result; by introducing the coupling of Eqs. (31)–(33), this problem is systematically addressed at the discrete level.

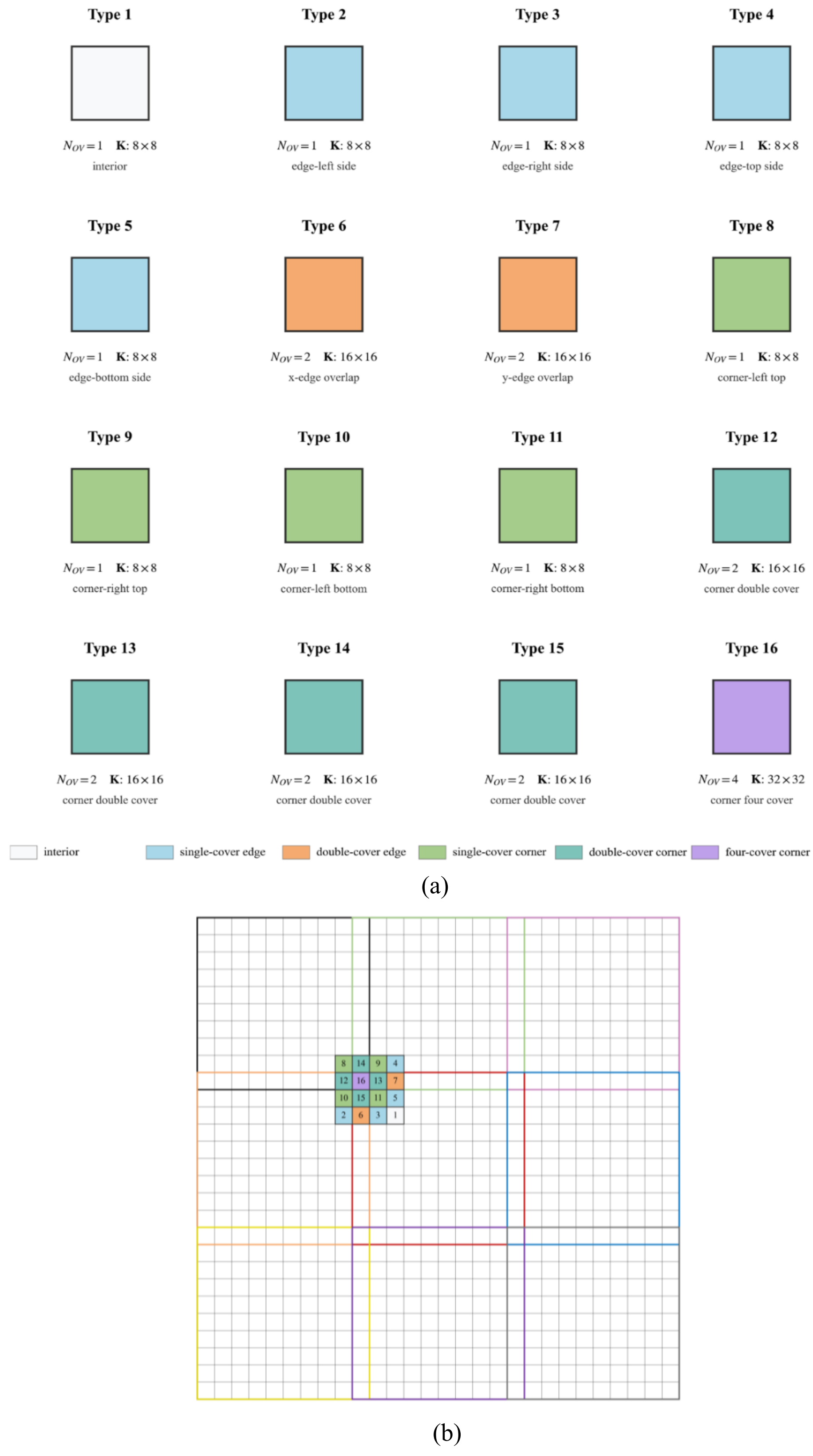

Fig.10. The 16 types of fine elements. (a) Catalogue of the 16 fine element types. (b) Locations of the 16 types in the overlapped Cartesian grid.

From the perspective of the computational workflow, the above coupling process can be summarized as: fine-scale compatible assembly, corner node-space projection and solution, and fine-scale displacement recovery. First, the overlapping finite elements complete the stiffness assembly on the unified Cartesian fine mesh, whose unknown form corresponds to the degrees of freedom of all fine-scale nodes; subsequently, through the oversampling mapping $\widetilde{\boldsymbol{\Phi}}_{\mathrm{global}}$, the fine-scale system is projected into the low-dimensional subspace spanned by the substructures' corner node degrees of freedom, thereby reducing the size of the unknowns solved online from the fine-mesh degrees of freedom down to the corner node degrees of freedom. In the two-dimensional case, each substructure retains only four corner nodes and eight degrees of freedom, so the size of the overall global unknowns is comparable to that of the classical four-corner node substructure method. After solving for the corner node displacements $\boldsymbol{u}_{\mathrm{c}}^{j}$, the internal fine-scale displacement of each substructure is recovered using Eq. (16), and combined into the final $\boldsymbol{u}_{\mathrm{global}}$ according to the partition-of-unity weighting functions. The recovery process is essentially a multiplication between a sparse matrix and a vector, and its computational cost is typically very low compared to solving the global system of linear equations. In other words, the overlapping finite elements ensure that the fine-scale space prior to projection is compatible, while the oversampling mapping ensures that the basis used for the projection is capable of capturing complex local deformations; neither one can be dispensed with the other.

This coupling approach has direct implications for mechanical analysis and topology optimization. For mechanical analysis, it enables the recovery of a continuous displacement field approaching that of fine-scale finite elements while retaining only corner node degrees of freedom, thereby improving the reliability of strain and strain energy computation near interfaces. For topology optimization, the material distribution changes frequently during the iteration process; if each substructure were to solve the oversampling local problem online, the computational cost would accumulate rapidly. However, once the local mapping can be rapidly provided by a data-driven model, overlapping finite element assembly and corner node reduced-order solution become the dominant parts of the online cost, and these two parts maintain high efficiency due to the Cartesian mesh and typed element stiffness precomputation. More importantly, the machine learning module predicts local oversampling shape functions, rather than the global solution of a specific boundary value problem. Therefore, the trained model can be reused under different structural configurations, loading forms, and boundary conditions, thus together with the overlapping finite element reduced-order framework forming a problem-independent acceleration path.

# 4 Fast Prediction of Oversampled Numerical basis functions Using U-Net

Section 2 provided the local construction method for oversampled numerical basis functions, while Section 3 established the global displacement compatibility mechanism through partition-of-unity-based overlapping finite elements and embedded the local oversampling mapping into this compatible space. Thus, the reduced-order analysis chain of PIML-OFEM has been established; in this computational framework, the main bottleneck for online efficiency no longer arises solely from the scale of the global degrees of freedom, but from the fact that each substructure needs to re-solve a local boundary value problem based on the material distribution in its oversampling region in order to obtain the mapping matrix from corner node degrees of freedom to fine-scale degrees of freedom. When the number of substructures is large or the material distribution changes frequently during the iteration process, constructing oversampled numerical basis functions substructure by substructure would significantly erode the efficiency advantage of the reduced-order model.

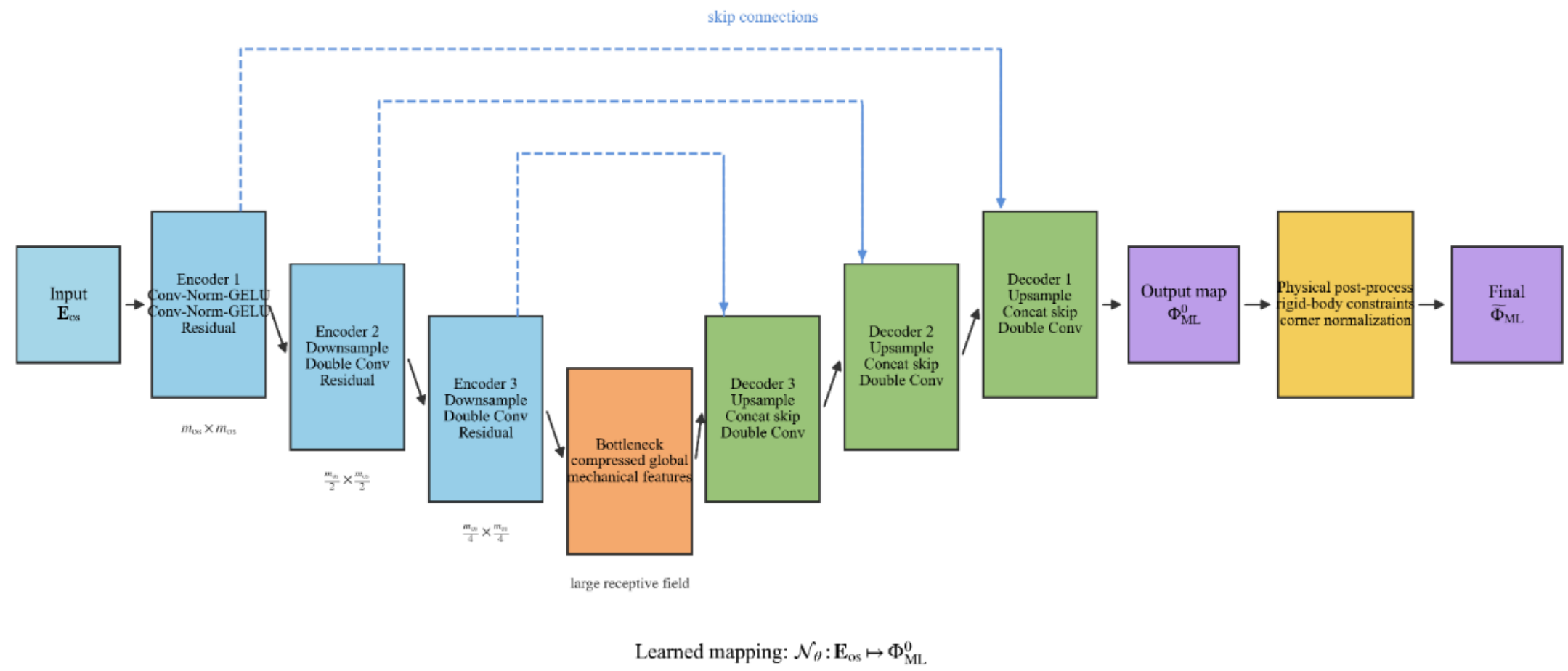

Fig. 11 U-Net network architecture.

To reduce this part of the computational cost, this paper constructs a U-Net-based data-driven prediction model to learn the deterministic mapping from the oversampling region material distribution to the local numerical basis functions. The network architecture is shown in Fig. 11. It should be emphasized that this model does not directly predict the global displacement field or the global stiffness matrix, but rather predicts the local oversampled numerical basis functions defined in Section 2; the predicted results, after rigid-body mode completion and corner node normalization, are embedded as local reduced-order mappings into the overlapping finite element compatibility framework of Section 3. Compared with purely fully connected networks, U-Net can exploit its convolutional structure to extract local spatial features from the Young's modulus field and, through the encoder-decoder architecture, simultaneously account for local details and larger-scale response characteristics, and is thus suitable for handling material distributions with voids, interfaces, and strong heterogeneity.

**4.1 U-Net Network Structure**

Let the number of fine grid elements in the x and y directions inside the substructure be denoted by m, and let the number of fine grid layers extended outward for oversampling be l. The number of fine grid elements in the x, y directions within the oversampling region is then:

$$m_{\mathrm{os}} = m + 2l. \tag{34}$$

For linear elastic problems, the network input is taken as the Young's modulus matrix $\boldsymbol{E}_{\mathrm{os}} \in \mathbb{R}^{m_{\mathrm{os}} \times m_{\mathrm{os}}}$ of the $m_{\mathrm{os}} \times m_{\mathrm{os}}$ fine-scale elements within the oversampling region. With the Poisson's ratio fixed, $\boldsymbol{E}_{\mathrm{os}}$ specified the local stiffness distribution that determines the value of the oversampled numerical basis functions. Since the different types of substructures in Section 3.1 have different local mesh topologies, the number of their fine-scale degrees of freedom $ns$ may also differ. Taking the "internal" type substructure as an example, its internal fine-scale degree-of-freedom count is:

$$ns = 2(m-1)^2. \tag{35}$$

At this point, the network output is organized as a partial numerical basis function matrix, in which each fine-scale node corresponds to two shape function components in the x and y directions, and the column direction corresponds to the combination of corner node degrees of freedom after removing the rigid-body mode constraints. This output is then rearranged into the degree-of-freedom format $\boldsymbol{\Phi}_{\mathrm{ML}}^{0} \in \mathbb{R}^{ns \times (nc-3)}$, and completed into the full local numerical basis function matrix $\boldsymbol{\Phi}_{\mathrm{ML}} \in \mathbb{R}^{ns \times nc}$ by supplementing it with the rigid-body mode constraints; the specific completion procedure can be found in [37]. This operation ensures that the resulting predicted shape functions satisfy the necessary rigid-body mode constraints and corner node normalization requirements through hard physical constraints, thereby guaranteeing that no spurious strain energy is generated when the substructure undergoes rigid-body translation or rotation.

Before entering the network, $\boldsymbol{E}_{\mathrm{os}}$ is treated as a single-channel two-dimensional feature map.

To facilitate multi-level downsampling and upsampling operations and to mitigate the information truncation effect at the convolution boundaries, appropriate padding can be applied at the input boundaries to obtain an expanded input feature map. Subsequently, the network gradually extracts spatial features induced by the material distribution through several encoder modules.

The encoder part consists of three levels of convolutional down-sampling modules. Each encoder level includes 2 convolutional layers, a normalization layer, and a GELU activation function, and residual connections are introduced to improve the training stability of deep networks. As the network depth increases, the spatial scale of the feature maps gradually decreases while the number of channels gradually increases, enabling the network to capture the influence of the material distribution on the overall local response at a smaller spatial scale. At the tail end of the encoder, a bottleneck layer is set up to further extract global features after compression. This layer is located at the lowest spatial resolution of the U-Net, has a relatively large receptive field, and can synthesize the coupling effects of material measures at different positions within the oversampling region on the local shape function. Denote the output feature of the $k$-th level encoder as $F_{enc}^{k}$, it contains a progressively abstracted representation, from the local material distribution to relatively higher-level mechanical response features.

The decoder part is roughly symmetric to the encoder. It progressively restores the spatial resolution of the feature maps through up-sampling operations, and at each level, the corresponding encoder-level features are passed to the decoder via skip connections. This skip connection can be expressed as:

$$F_{dec}^{k} = D^{k}\left(F_{dec}^{k+1}, F_{enc}^{k}\right), \tag{36}$$

where $D^{(k)}(\cdot)$ denotes the $k$-th level decoder operation. In this way, the network can utilize both the global response information in the deep features and preserve the local spatial details in the shallow features, thereby improving its predictive capability for numerical basis functions under complex heterogeneous material distributions.

After the decoder, an output mapping layer is set up to convert the restored feature map into the predicted result $\boldsymbol{\Phi}_{\mathrm{ML}}^{0}$ in degree-of-freedom format. Since this matrix does not yet fully satisfy the rigid-body mode and corner node interpolation constraints, this paper subsequently performs the physical constraint completion and the corner node normalization operation from Section 2.2 after the network output, obtaining the final $\boldsymbol{\Phi}_{\mathrm{ML}}$ used for substructure condensation.

In summary, the U-Net network adopted in this paper essentially constructs the following mapping relationship:

$$N_{\boldsymbol{\theta}}: \boldsymbol{E}_{\mathrm{os}} \longmapsto \boldsymbol{\Phi}_{\mathrm{ML}}^{0}, \tag{37}$$

where $\mathbf{N}_{\boldsymbol{\theta}}$ denotes the U-Net network with parameters $\boldsymbol{\theta}$. This mapping directly converts the material distribution within the local oversampling region into the corresponding numerical basis function, thereby avoiding the need to repeatedly solve the local oversampling boundary value problem for every substructure during the online stage. Because the local boundary-value problem defines a deterministic mapping from $\boldsymbol{E}_{\mathrm{os}}$ to $\boldsymbol{\Phi}_{\mathrm{ML}}^{0}$, this mapping is particularly suitable for approximation using a machine learning model. Because the convolutional structure preserves the spatial neighborhood information of the material distribution, and the skip connections take into account both local details and global features, this network can approximate, at relatively low computational cost, the local numerical basis functions obtained by traditional oversampling methods. Samples are generated respectively for nine different types of substructures to conduct training, yielding nine machine learning models.

**4.2 Training Data and Loss Function**

The training data is generated offline according to the traditional oversampling computation workflow described in Section 2.2. For each material sample, the local boundary value problem is first solved within the corresponding oversampling region, and then the target numerical basis function is extracted and normalized. For substructure learning with m=10, the dataset contains a total of 20,000 samples, of which 19,000 samples are used for training and 1,000 samples for validation. Each sample includes the oversampling elastic modulus input matrix $\boldsymbol{E}_{\mathrm{os}}$, the corresponding partial numerical basis function label $\boldsymbol{\Phi}_{\mathrm{ML}}^{0}$, and the local stiffness matrix $\mathbf{K}_{\mathrm{loc}}$ used for stiffness-fidelity constraints. For each sample, an independent random density field is drawn as $E_{\mathrm{sample}} \sim \mathcal{U}[10^{-6}, 1]$. The random seed is fixed at 42 for reproducibility. Both the input and output are standardized based on statistics from the training set before training, and the same standardization parameters are reused during validation and online prediction.

This paper adopts a composite loss function jointly composed of the data error and the stiffness-fidelity error:

$$L = \lambda_{\mathrm{data}} \, L_{\mathrm{data}} + \lambda_{\mathrm{K}} \, L_{\mathrm{K}}, \tag{38}$$

where $L_{\mathrm{data}}$ uses the Smooth L1 loss function, used to measure the per-degree-of-freedom error between the predicted shape function and the offline oversampling label, and to improve the model's robustness to locally anomalous samples; $L_{\mathrm{K}}$ is the stiffness-fidelity term, used to measure the relative error between the reduced-order stiffness matrix induced by the predicted shape function and the true reduced-order stiffness matrix:

$$L_{\mathrm{K}} = \frac{\left\| (\widetilde{\boldsymbol{\Phi}}^{j}{}_{\mathrm{ML}})^{\mathsf{T}} \mathbf{K}_{\mathrm{loc}}^{j} \widetilde{\boldsymbol{\Phi}}^{j}{}_{\mathrm{ML}} - \widetilde{\boldsymbol{\Phi}}^{j^{T}} \mathbf{K}_{\mathrm{loc}}^{j} \widetilde{\boldsymbol{\Phi}}^{j} \right\|_{\mathrm{F}}}{\left\| \widetilde{\boldsymbol{\Phi}}^{j^{T}} \mathbf{K}_{\mathrm{loc}}^{j} \widetilde{\boldsymbol{\Phi}}^{j} \right\|_{\mathrm{F}} + \varepsilon}. \tag{39}$$

where $\varepsilon = 10^{-8}$ is used to avoid numerical instability caused by an excessively small denominator, and $\left\| . \right\|_{\mathrm{F}}$ denotes the Frobenius norm of the matrix. During training, $\lambda_{\mathrm{data}} = 1.0, \lambda_{\mathrm{K}} = 0.1$ are used, and only the data term is used for warm-up during the first 50 epochs, after which the stiffness constraint term is enabled. This strategy avoids the stiffness error exerting excessive interference on the optimization direction during the early training stage, when the predicted shape functions are still unstable.

The optimizer used is AdamW, with an initial learning rate of $2 \times 10^{-4}$, weight decay of $10^{-5}$, and a batch size of 64. The learning rate adopts a cosine annealing restart strategy, with a minimum learning rate of $10^{-6}$. During training, the validation mean absolute error (validation MAE) is used as the model selection criterion, and an early stopping mechanism is set up.

**4.3 Training Results**

The model was trained on an NVIDIA GeForce RTX 4090 GPU. The model converged rapidly in the early stage of training: the validation MAE at the 1st epoch was 0.089689, dropping to 0.021373 at the 100th epoch; the error then further decreased slowly and stabilized in the later stages. Training triggered the early stopping mechanism at the 897th epoch, after which the model weights with the best validation set metric were restored. The final training loss was 0.009496, and the final validation loss was 0.006241; the best validation loss achieved during training was 0.006126, with a best validation MAE of 0.016082.

The training and validation histories are shown in Fig. 12. Their convergence indicates that the U-Net learns the nonlinear mapping from the local modulus field to the partial oversampled basis with stable validation performance. The stiffness-fidelity term also remains small during the later stages of training, showing that the predicted bases reproduce not only the reference entries but also the associated local reduced stiffness matrices.

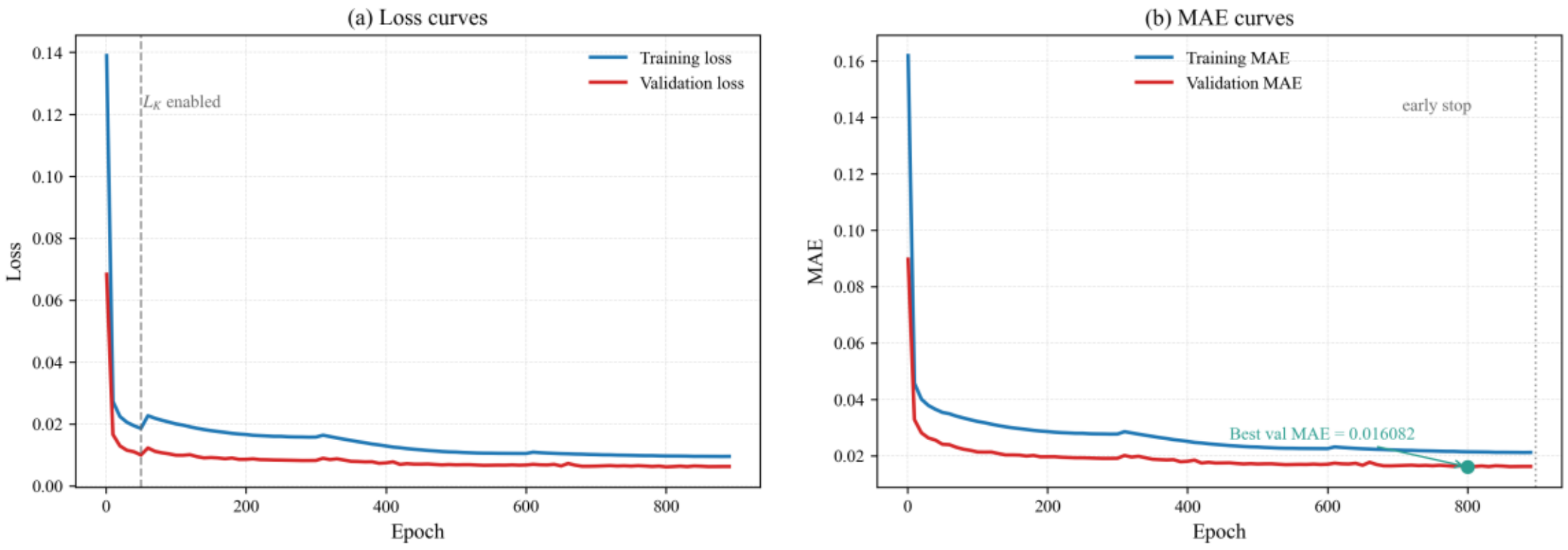

Fig. 12 U-Net model training and validation curves.

### 4.4 Online Prediction and Computational Workflow

In the online condensation-solution stage, the program first extracts the corresponding material field for each substructure according to the oversampling region definition of Section 2; for substructures near the overall boundary, the virtual material extension strategy from Section 2.2 is still used to supplement the outer extended region. Subsequently, the inputs of each substructure are organized into batch data by type and size and fed into the U-Net model. For substructures with approximately uniform Young's modulus —if the difference between the maximum and minimum Young's modulus within the oversampling region is less than $10^{-3}$, the pre-computed uniform-modulus shape function is directly reused, thereby further reducing the number of neural network calls; for heterogeneous substructures, the U-Net is used for batch prediction to obtain their local numerical basis function matrices.

After the prediction is complete, rigid-body mode completion and corner node normalization are performed separately for each substructure, yielding $\mathbf{\Phi}_{\mathrm{ML}}^{j}$, which can be used for condensation analysis. Subsequently, the local mappings of all substructures are assembled into a global reduced-order mapping matrix, and combined with the global compatible stiffness matrix formed by the partition-of-unity overlapping finite elements in Section 3, the condensed stiffness matrix $\mathbf{K}_{\mathrm{global}}^{\mathrm{c}}$ is constructed according to Eq. (33). Finally, the condensed system is solved, and the predicted shape functions are used to recover the fine-scale nodal displacements from the solved corner-node degrees of freedom. Therefore, the U-Net module only changes the way the oversampled numerical basis functions are obtained, mainly replacing the most time-consuming step of solving the local oversampling finite element boundary value problem; it does not change the displacement compatibility mechanism of the overlapping finite elements, nor does it change the overall finite element assembly, boundary condition application, and global equilibrium equation solving procedures.

## 5 Numerical Examples

This section evaluates the analysis accuracy, computational efficiency, and coupling capability with topology optimization of the proposed PIML-OFEM method through three numerical examples. All examples use the classical fine-scale finite element method (FEM) as the reference solution and are compared under the same material distributions, identical boundary conditions, and consistent post-processing conventions. The first example adopts a cantilever beam problem to validate the mesh coordination of the overlapping substructure discretization format, the local deformation recovery capability without presupposing substructure boundary displacement patterns, and the generalization capability of the trained machine learning model under different material distributions. The second example considers heterogeneous structures containing complex void topologies, used to further test the method's recovery accuracy for displacement-derived quantities and element strain energy—higher-order responses—in regions with complex void boundaries and local stress concentrations. The third example embeds PIML-OFEM into the SIMP topology optimization framework to examine its stability, efficiency, and capacity to preserve fine-scale topological features across multiple rounds of iterative analysis. The three examples respectively correspond to three key issues in the chain of the proposed method: whether the interface compatibility is reliable, whether local higher-order responses can be accurately recovered, and whether the analysis kernel can support high-resolution optimization iterations. In this section, unless otherwise stated, the plane stress assumption is adopted throughout, with Poisson’s ratio $\nu = 0.3$, and a dimensionless convention is used.

All computations were performed on a workstation equipped with a 13th-generation Intel Core i9-13900K CPU, an NVIDIA GeForce RTX 4090 GPU, and 128 GB of RAM. The same trained model was used in all examples. It was trained for two-dimensional substructures containing $10 \times 10$ fine-scale elements, with the oversampling domain extended by six fine-scale element layers on each side. During online analysis, the precomputed uniform-material basis was reused for oversampling domains with approximately uniform material properties, whereas heterogeneous domains were evaluated in batches by the U-Net. No model was retrained for a particular structure, boundary condition, or load; the examples therefore test transfer across these global problem settings within the local discretization and material range used for training.

To further improve the local accuracy in higher-order response evaluation and topology optimization sensitivity computation, after completing the condensed system solution and fine-scale displacement reconstruction, this paper uses the reconstructed displacement as the initial field and performs a small number of residual corrections on the original fine-scale equilibrium equations. This step does not alter the previously constructed oversampled numerical basis functions or the overlapping finite element compatibility space, but rather exploits the high-quality initial displacement field already provided by PIML-OFEM to rapidly reduce the equilibrium residual. Specifically, the preconditioned conjugate gradient method (PCG) is adopted with the diagonal terms of the global stiffness matrix as the Jacobi preconditioner. For problems with millions of two-dimensional elements, this correction process typically introduces only about 1-2 s of additional computation time, yet it can effectively improve the post-processing quality of derived quantities such as strain, strain energy, and topology optimization sensitivity.

**5.1 Error Metrics**

To quantitatively evaluate the consistency between the PIML-OFEM and classical FEM solutions, this paper separately computes the nodal displacement error and the elemental strain energy error. Let the nodal displacement vector obtained by PIML-OFEM be $\boldsymbol{u}_{\text{est}}$, and the classical FEM reference solution be $\boldsymbol{u}_{\text{ref}}$. For the $\beta$-th node, the predicted displacement can be written as $\boldsymbol{u}_{\text{est}}^{\beta} = \left(u_{\text{est}}^{\beta}, v_{\text{est}}^{\beta}\right)^{\mathrm{T}}$, and the reference displacement as $\boldsymbol{u}_{\text{ref}}^{\beta}$. The nodal relative displacement error is defined as $\eta_{\mathrm{u}}^{\beta} = \frac{\left\|\boldsymbol{u}_{\text{est}}^{\beta} - \boldsymbol{u}_{\text{ref}}^{\beta}\right\|_2}{\left\|\boldsymbol{u}_{\text{ref}}^{\beta}\right\|_2 + \varepsilon_{\mathrm{u}}}$, where $\varepsilon_{\mathrm{u}}$ represents a small scalar, set to $10^{-12}$ in the present study to avoid singularity in the relative error where the reference displacement approaches zero. The domain-averaged nodal relative displacement error is defined as $\bar{\eta}_{\mathrm{u}} = \frac{\sum_{\beta=1}^{N_{\mathrm{n}}} \eta_{\mathrm{u}}^{\beta}}{N_{\mathrm{n}}}$, where $N_{\mathrm{n}}$ is the total number of nodes. This metric reflects the deviation of the overall displacement field relative to the reference solution, in an average sense.

In addition to the displacement error, this paper also uses the elemental strain energy error to evaluate the recovery accuracy of displacement derivative quantities. For the $e$-th element, the standard Q4 element stiffness matrix $\mathbf{K}_{e}^{\mathrm{Q4}}$ (evaluated at unit Young's modulus) is uniformly used for post-processing, and the strain energy of the $e$-th element can be defined as $SE_e =$ $\frac{1}{2}\, E_e\, (\boldsymbol{u}_e)^\top\, \mathbf{K}_e^{\mathrm{Q4}}\, \boldsymbol{u}_e$, where $E_e$ is the relative material Young's modulus of the $e$-th element, and $\boldsymbol{u}_e$ is the corresponding elemental degree-of-freedom displacement vector. It should be emphasized that, in computing the strain energy error, different methods all use the same standard Q4 stiffness matrix for post-processing, so the error mainly reflects the difference in the displacement field itself, rather than differences in discretization scheme or element stiffness definition. The relative elemental strain energy error is defined as $\eta_{\mathrm{SE}}^e = \frac{|SE_{\mathrm{est}}^e - SE_{\mathrm{ref}}^e|}{|SE_{\mathrm{ref}}^e| + \varepsilon_{\mathrm{SE}}}$, where $\varepsilon_{\mathrm{SE}} = 10^{-12} \max_e |SE_{\mathrm{ref}}^e|$. This metric is particularly critical for topology optimization, since in compliance minimization problems the sensitivity is directly related to the elemental strain energy, and systematic strain energy errors may accumulate during iteration and affect the final configuration.

**5.2 Cantilever Beam Example**

Consider the two-dimensional cantilever beam of dimensions 2×1 shown in Fig. 13. The left edge is fixed, and a downward concentrated force $\boldsymbol{f} = 1$ is applied at the lower-right corner. The concentrated load generates steep local displacement gradients near its point of application. This example is therefore suitable for evaluating whether a corner-node-based reduced model can recover a nonuniform boundary-displacement profile and localized deformation without introducing additional boundary control degrees of freedom.

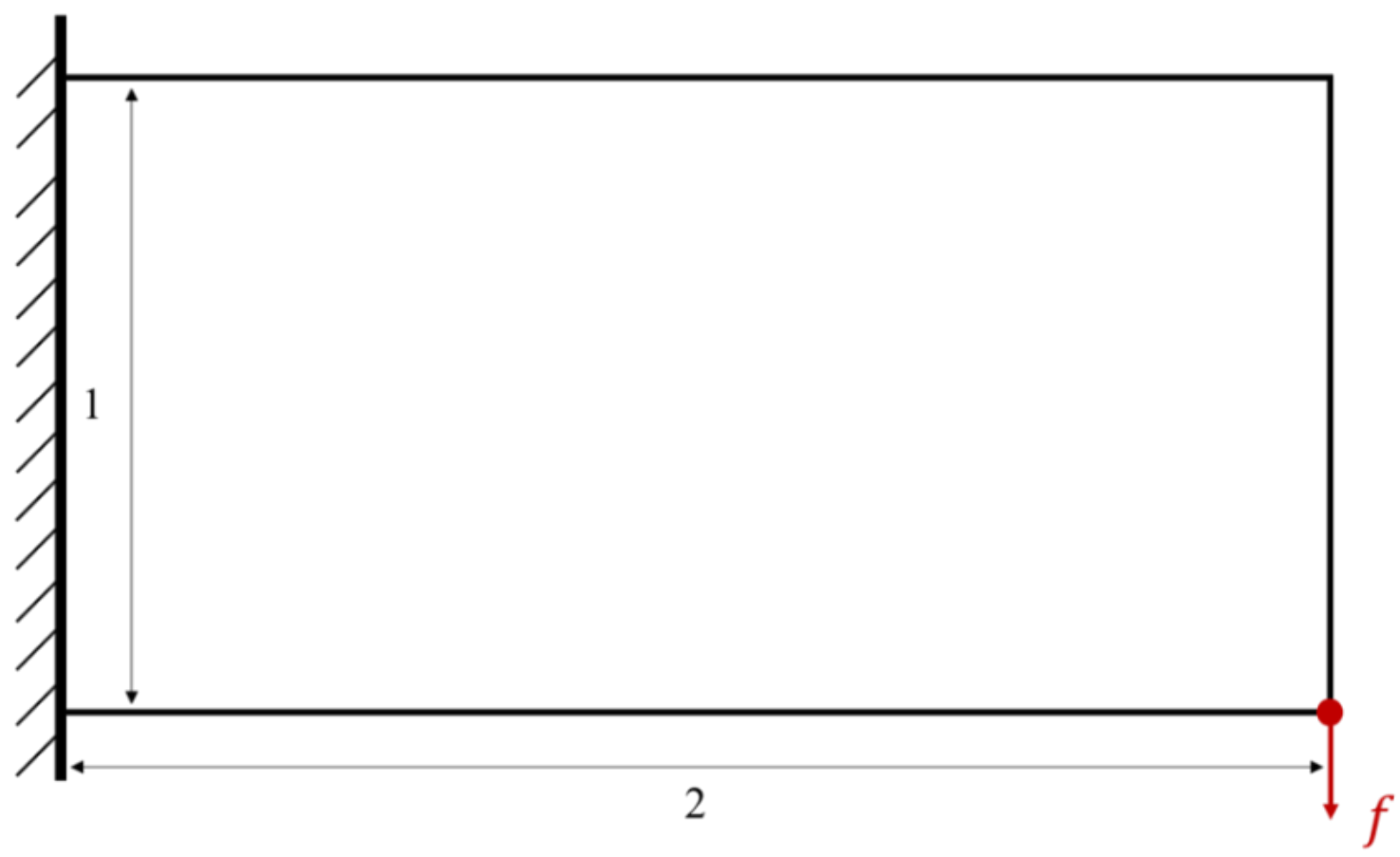

Fig.13. Boundary conditions of the cantilever beam example.

First consider a fully solid structure, i.e., the entire design domain is filled with material of Young's modulus $E_{\max} = 1$. To verify the effectiveness of the proposed overlapping substructure discretization scheme, analyses are performed using both the classical FEM mesh and the PIML-OFEM mesh. PIML-OFEM divides the structure into 8×4 mutually overlapping substructures, each containing 10×10 fine-scale elements; the equivalent fine-scale resolution obtained after overlapping is consistent with classical FEM, both being 73×37. As shown in Fig. 14, the deformed mesh obtained by PIML-OFEM is highly consistent with the classical FEM result, and no visible displacement jump or mesh distortion was observed at substructure interfaces. This result indicates that the partition-of-unity-based overlapping compatibility mechanism proposed in Section 3 can effectively eliminate the interface incompatibility problems commonly seen in traditional non-overlapping substructure assembly.

It is worth noting that, in this example, the displacement field near the load application point has strong locality and sharp-corner features. Traditional PIML or substructure methods based on the boundary linear assumption typically find it difficult to accurately reconstruct such localized deformation while retaining only corner node degrees of freedom; if the boundary interpolation order is increased to improve accuracy, this inevitably increases the number of retained degrees of freedom and reduces the sparsity of the condensed stiffness matrix. In contrast, PIML-OFEM still retains only four corner nodes—eight degrees of freedom in total—at the level of each substructure, but through the combination of oversampled numerical basis functions and overlapping compatible interpolation, it is able to recover the complex displacement pattern near the load point without introducing additional boundary degrees of freedom.

Subsequently, under the same cantilever beam boundary conditions, heterogeneous material distributions are further considered. Fig. 15(a) presents a double-circular-hole example, in which the hole region uses low-stiffness material $E_{\mathrm{min}} = 10^{-6}$, and the remaining region is solid material $E_{\mathrm{max}} = 1$. The computational domain is discretized with 1801×901 fine-scale elements. The classical FEM analysis takes 12.39 s; the PIML-OFEM analysis takes 2.64 s, of which solving the condensed linear system takes only 1.04 s. Using classical FEM as the reference solution, the average nodal relative displacement error of PIML-OFEM is $\bar{\eta}_{\mathrm{u}} = 0.0140$. At the location of maximum response, the reference solution gives a vertical displacement of $v_{\mathrm{ref}} = 71.30$, while PIML-OFEM predicts $v_{\mathrm{est}}^{\beta} = 70.97$, showing good agreement. The displacement contours and error contours in Fig. 15(a) further show that this method can accurately capture the local perturbation induced by the holes and the overall bending response.

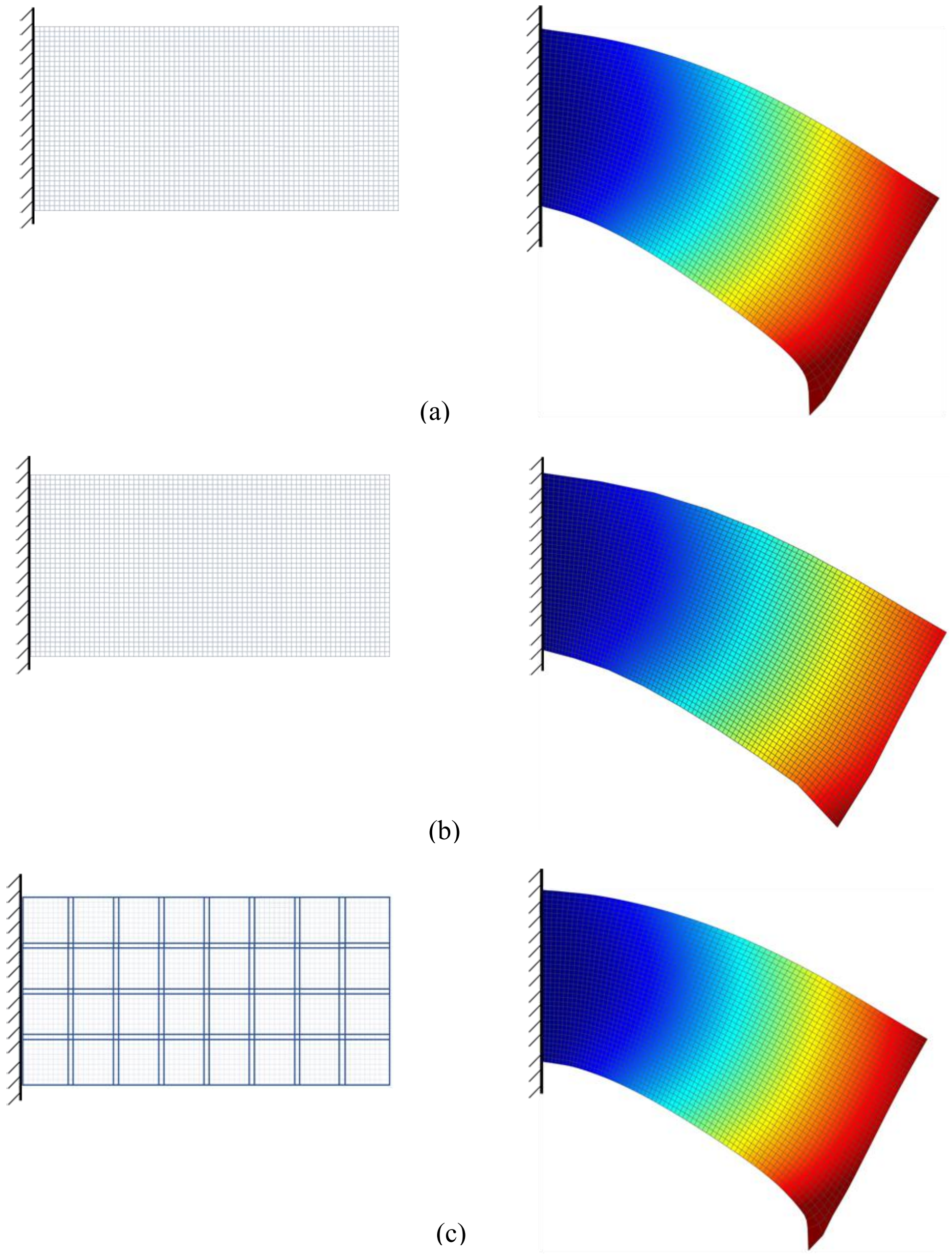

Fig. 14. Comparison of cantilever beam deformation under different discretization and reduction strategies. (a) Result of classical FEM analysis. (b) Result of the PIML method based on the linear boundary displacement assumption for substructures [27-28]. (c) Result of PIML-OFEM analysis, in which the overall structure is discretized into 8×4 mutually overlapping substructures, each containing 10×10 fine-scale elements; the equivalent fine-scale resolution after overlapping is consistent with panel (a), both being 73×37.

To more rigorously examine the robustness of the method under strongly heterogeneous conditions, a random material modulus example is further constructed. In this example, taking the fine-scale element as the smallest unit, a random Young's modulus in the range $[E_{\min}, E_{\max}]$ is assigned throughout the domain, forming a highly oscillatory material field with no obvious spatial correlation. This setting imposes higher demands on the prediction accuracy of the local oversampled shape functions and the internal displacement recovery capability of the substructures. The computational domain is likewise discretized with 1801×901 fine-scale elements. Classical FEM analysis takes 13.80 s; PIML-OFEM analysis takes 5.83 s, of which solving the condensed linear system takes only 0.976 s. In this example, the average nodal relative displacement error is $\bar{\eta}_u = 0.0692$. At the location of maximum response, the classical FEM solution is $v_{\mathrm{ref}} = 256.29$, while the PIML-OFEM result is $v_{\mathrm{est}}^{\beta} = 248.30$.

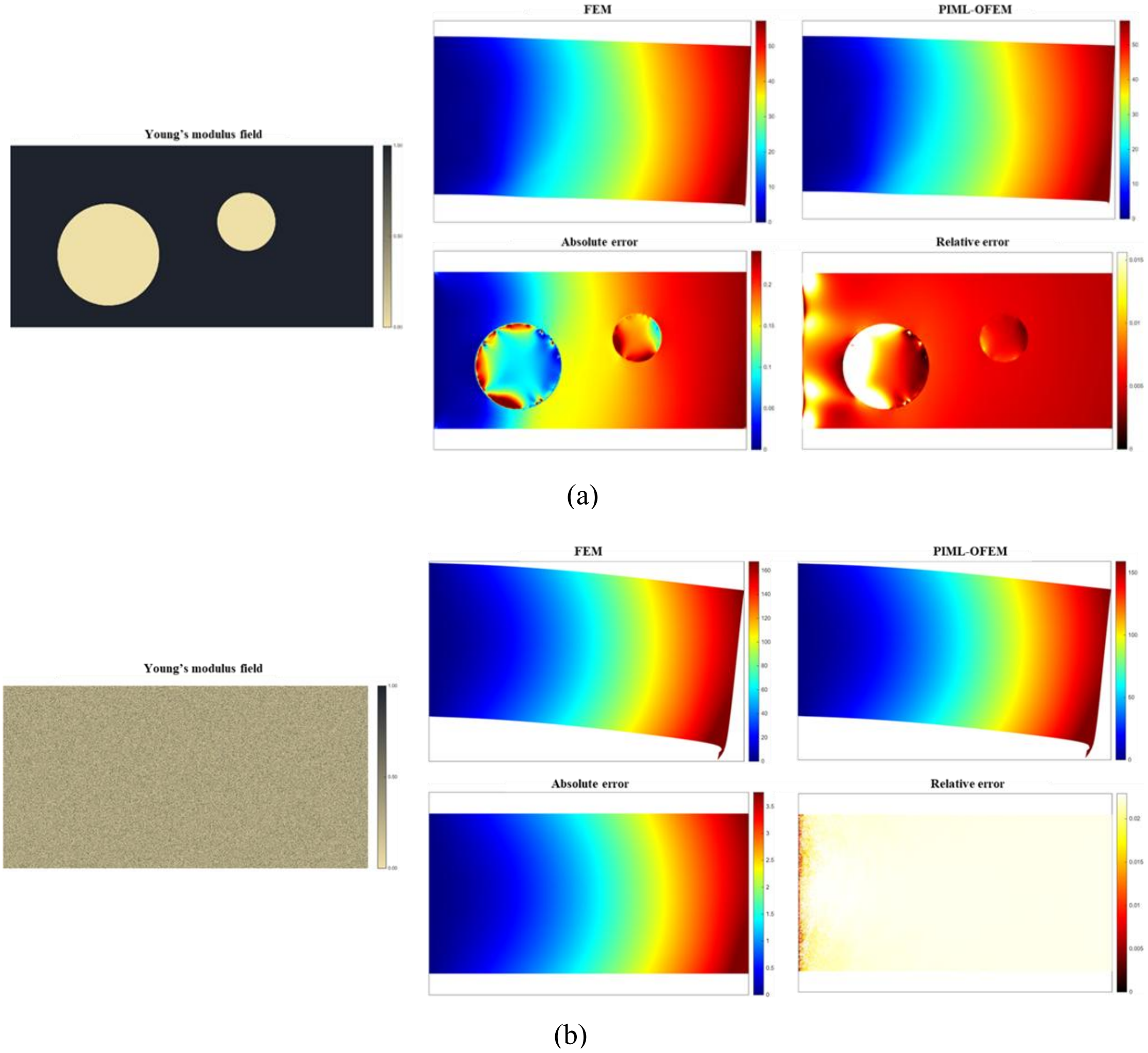

Fig.15. Displacement prediction results and error distributions in the heterogeneous cantilever beam example. (a) Double-circular-hole example; (b) random Young's modulus example. Each group of results respectively shows the deformation and displacement contours of classical FEM and PIML-OFEM, as well as the absolute and relative error contours between the two.

Compared with the double-circular-hole example, in the random Young's modulus example almost all internal substructures exhibit heterogeneous material distributions, so the uniform substructure shape function reuse mechanism cannot be extensively utilized, resulting in increased PIML-OFEM online time. Nevertheless, even under such strongly heterogeneous and fine-scale randomly oscillating material fields, PIML-OFEM can still maintain acceptable overall displacement accuracy while significantly shortening the linear system solution time. The deformation contours in Fig. 15(b) show that the complex local deformation induced by the random material distribution near the load application point can be well recovered, indicating that the trained U-Net model has learned not only the smooth response under uniform or weakly heterogeneous materials, but also possesses the capability to handle strong local material perturbations.

**5.3 Example with Complex Void Topology**

The second example considers the complex topology structure shown in Fig. 16. The computational domain size is 4×1, containing internally several irregular holes, including a tilted elliptical hole, circular holes, a diamond-shaped hole, and randomly distributed elliptical holes. The left end of the structure is fixed, and a vertical concentrated force of magnitude $\boldsymbol{f} = 1$ is applied at the bottom-right corner. The geometric boundary in this example is more complex, and there is noticeable local stress concentration near the holes, making it suitable for examining the reliability of PIML-OFEM in recovering complex topology and higher-order responses.

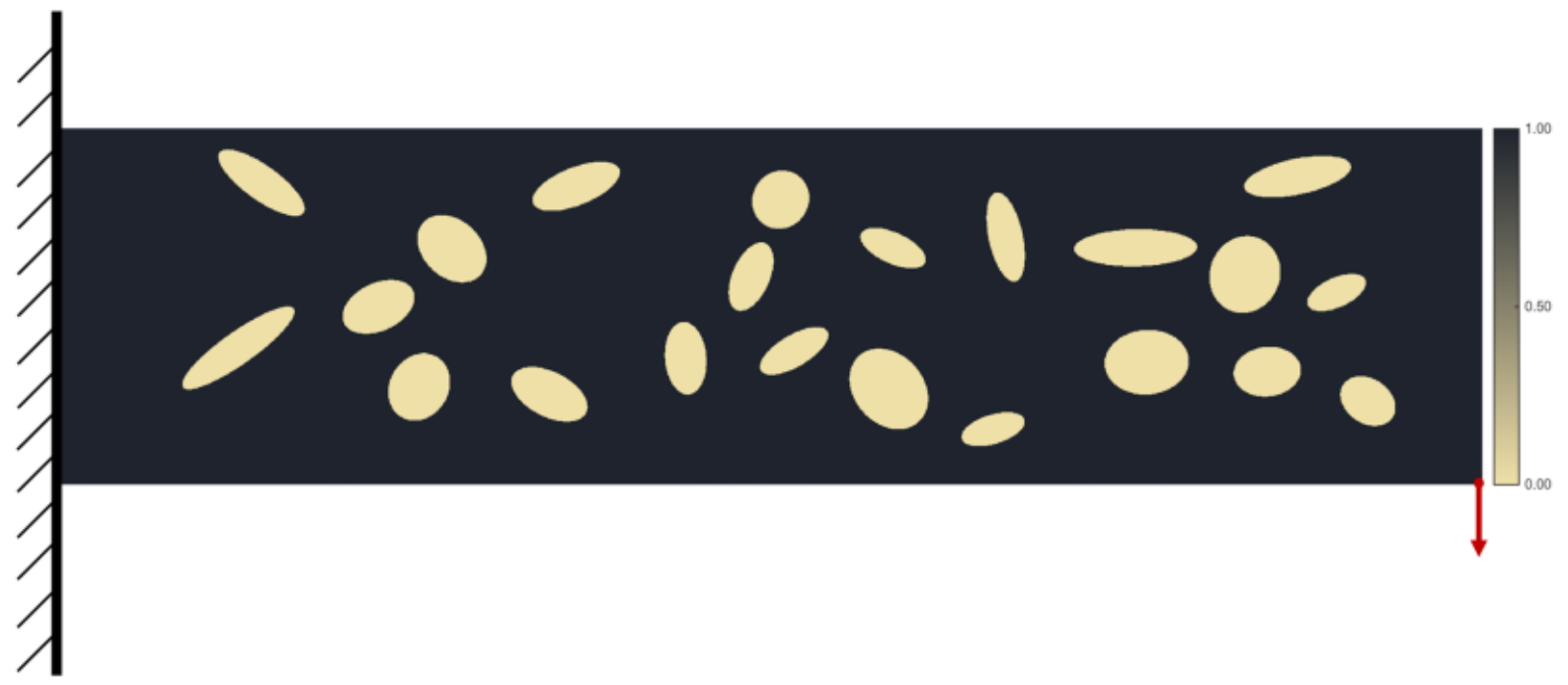

Fig.16. Material distribution and boundary conditions of the complex topology structure.

The overall structure is discretized with 3601×901 fine-scale elements. Classical FEM analysis takes 39.85 s; PIML-OFEM analysis takes 12.06 s, of which solving the condensed linear system takes 2.60 s. Using classical FEM as the reference solution, the average nodal relative displacement error of PIML-OFEM is $\bar{\eta}_{\mathrm{u}} = 0.0066$. At the location of maximum response, the reference solution gives a vertical displacement of $v_{\mathrm{ref}} = 344.70$, while the PIML-OFEM result is $v_{\mathrm{est}}^{\beta} = 343.00$. The displacement contours and error contours given in Fig. 17(a-d) show that the overall deformation patterns obtained by the two methods are highly consistent; the absolute displacement error is mainly concentrated near the load point and hole boundaries, while remaining at a low level in most of the solid regions. The relative displacement error is overall controlled within a small range, indicating that PIML-OFEM can accurately describe the global compliance and local deformation response in complex hole structures.

Table 5.1 further compares the computational efficiency and displacement prediction accuracy of different mechanical analysis methods in this complex topology example. As the reference solution, classical FEM has an analysis time of 39.85 s. The PIML-EMS-LBC method based on the linear boundary-displacement assumption reduces the analysis time to 11.14 s, but because its substructure boundary displacement is still constrained by the linear assumption, it is difficult to fully characterize nonlinear deformation features near complex hole boundaries and local stress concentration regions, giving an average displacement error of 0.0199. In contrast, the proposed PIML-OFEM method reduces the average displacement error to 0.0066 while only slightly increasing the analysis time to 12.06 s, indicating that the proposed oversampled numerical basis functions and overlapping compatibility mechanism can significantly improve displacement recovery accuracy in complex topology structures. On the other hand, the OFEM-SUB method (which, compared to PIML-OFEM, does not use U-Net to predict the oversampled numerical basis functions but instead uses an online computation approach) does not employ machine learning acceleration, and can achieve a lower average displacement error of 0.0046, but since solving the oversampled numerical basis functions requires solving the complete boundary value problem, its analysis time reaches as high as 146.77 s, far higher than the other methods. The above results indicate that PIML-OFEM significantly reduces online computational cost while approaching the accuracy of exact substructure analysis, thereby achieving a superior accuracy-efficiency balance in complex heterogeneous structure analysis.

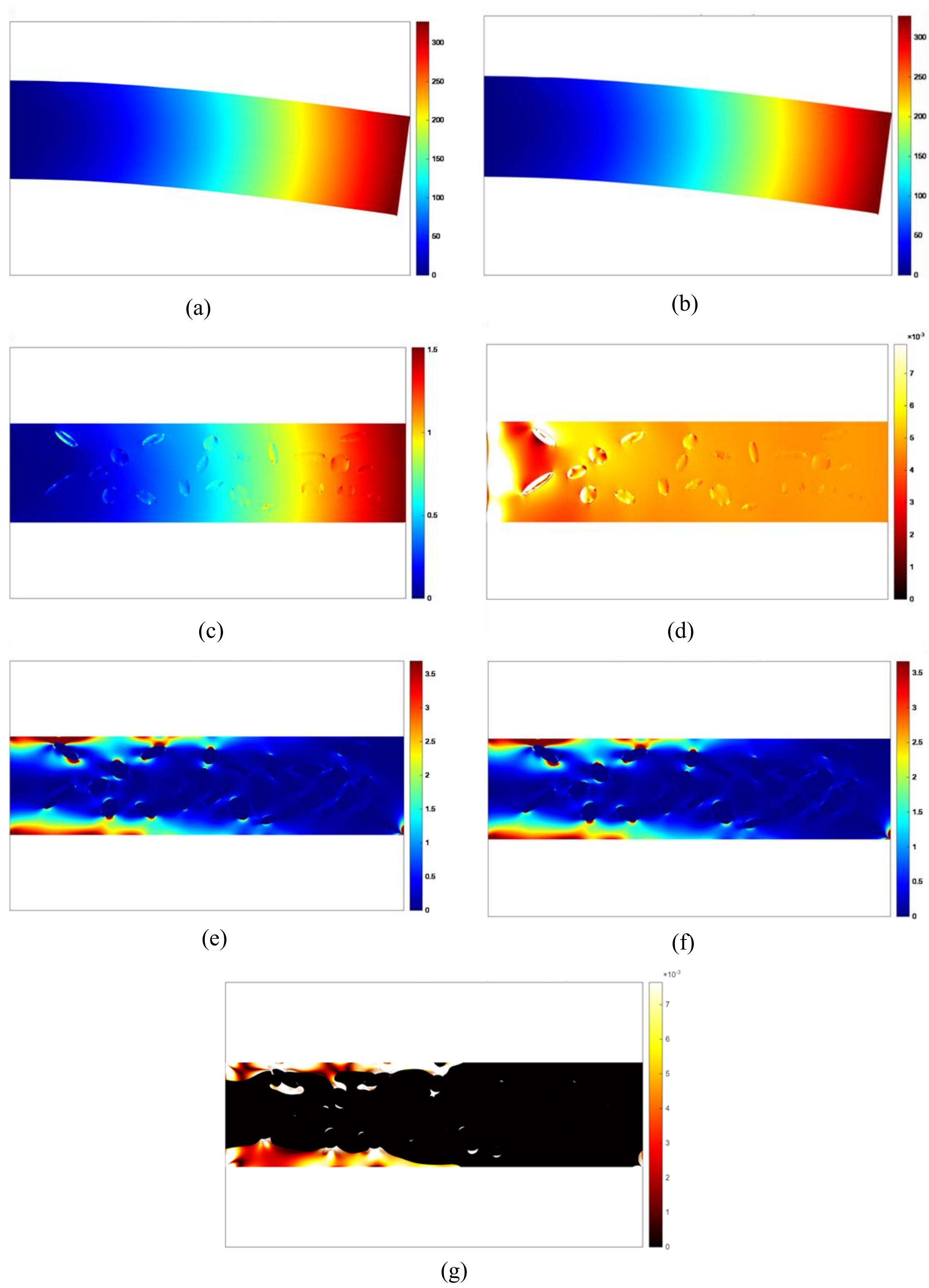

Fig.17. Mechanical analysis results of the complex topology structure under a concentrated force. (a) Deformation and displacement field contours obtained by classical FEM analysis. (b) Deformation and displacement field contours obtained by the PIML-OFEM method. (c) Absolute displacement error contour. (d) Relative displacement error contour. (e) Elemental strain energy

contour obtained by classical FEM. (f) Elemental strain energy contour obtained by the PIML-OFEM method. (g) Elemental strain energy relative error contour.

For topology optimization and structural strength assessment, comparing the displacement field alone is not sufficient. Since strain and strain energy involve spatial derivatives of the displacement field, they are more sensitive to local errors and interface incompatibility, and are therefore important indicators for evaluating the accuracy of reduced-order methods. Fig. 17(e-g) compare the elemental strain energy distribution and relative error obtained by classical FEM and PIML-OFEM. It can be seen that PIML-OFEM can not only recover the overall strain energy distribution, but also well capture the local high-energy regions near hole boundaries and load transfer paths. The relative error of elemental strain energy of the solid area is generally kept at a low level (below $8\times10^{-3}$), indicating that the proposed oversampled shape functions, overlapping compatibility mechanism, and small amount of residual correction can jointly improve the post-processing accuracy of displacement derivative quantities.

**Table 5.1** Comparison of computational efficiency and domain-averaged displacement error $\bar{\eta}_{\mathrm{u}}$ of different mechanical analysis methods in the complex topology example

| | FEM | PIML-EMS-LBC | PIML-OFEM | OFEM-SUB |
|---|---|---|---|---|
| Time (s) | 39.85 | 11.14 | 12.06 | 146.77 |
| Relative Error ( $\bar{\eta}_u$ ) | - | 0.0199 | 0.0066 | 0.0046 |

This example shows that the benefit of PIML-OFEM extends beyond the global displacement response. By enforcing compatibility across substructure interfaces and avoiding prescribed linear interpolation along target-substructure boundaries, the method reduces the local strain-energy bias associated with displacement jumps and restricted boundary kinematics. This improvement is important in topology optimization because element-wise strain energy directly enters the compliance sensitivity and persistent local errors can accumulate over successive design updates.

**5.4 Topology Optimization Example**

After verifying that PIML-OFEM has high accuracy for both displacement and strain energy, it is further embedded into a density-based topology optimization framework. This paper adopts the classical SIMP method for compliance minimization, replacing only the finite element analysis module with the proposed PIML-OFEM; the optimization variable update, sensitivity filtering, and volume constraint handling all follow the mature SIMP workflow. Therefore, the results in this section mainly reflect the accuracy, stability, and computational efficiency of the proposed mechanical analysis framework during topology optimization iteration, rather than the effect of introducing a new optimization algorithm.

First consider the Michell beam problem shown in Fig. 18(a), with a material volume fraction constraint of $\bar{v} = 0.5$, and the optimization objective is to minimize structural compliance. Fig. 18(b) shows the high-resolution topology result obtained in the literature based on the homogenization method under a 1600×800 background mesh [49]. Existing research has shown that, for this type of single-load-case problem, the local appearance of a rank-2-like microstructure is theoretically one of the reasonable optimal features.

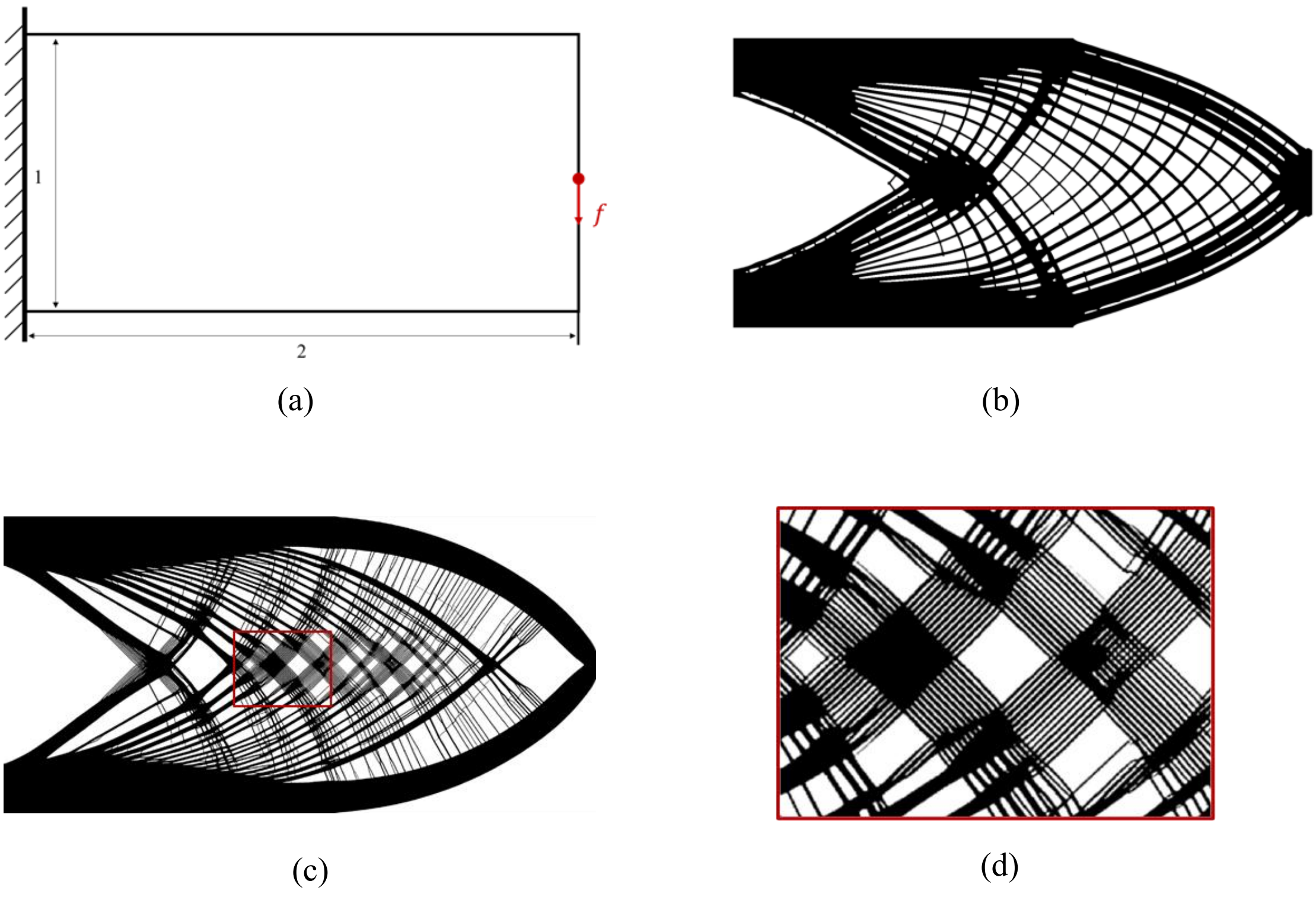

Fig.18. Michell beam topology optimization problem. (a) Design domain and boundary conditions. (b) Topology optimization result obtained using the de-homogenization topology optimization method at a resolution of 1600×800 [49]. (c) Topology optimization result obtained by directly optimizing with the proposed PIML-OFEM method in place of the finite element analysis in the classical SIMP algorithm. (d) Topological details of the structure shown in panel (c).

To examine the applicability of PIML-OFEM in high-resolution topology optimization, this paper discretizes the design domain into 3061×1531 fine-scale elements, and uses PIML-OFEM in place of traditional FEM to complete each equilibrium equation solution. The resulting optimization result is shown in Fig. 18(c), with the local zoom-in view shown in Fig. 18(d). It can be seen that the final structure forms rich fine-scale topological details, exhibiting, in some regions, multiscale material distribution exhibits local patterns that qualitatively resemble rank-2 microstructures. This result indicates that PIML-OFEM can support million-scale or even higher-resolution topology optimization analysis without significantly sacrificing analysis efficiency.

It should be noted that, in this example, the filter radius is only about √3 times the fine element size. For approaches that embed a reduced-order model or other approximate mechanical analysis model into topology optimization, if the local displacement field or elemental strain energy prediction has a large error, the compliance sensitivity will also show a bias in the corresponding region; such biases may manifest as element-scale or substructure-scale local oscillations during the iterative update. A larger filter radius can suppress the impact of such local errors on design variable updates by smoothing the density field or sensitivity field, and is therefore often used as a stabilization measure. However, the filter radius directly controls the minimum feature size of the optimized structure. Although an excessively large filter radius helps improve numerical stability, it weakens the fine detail expression capability that high-resolution optimization was originally intended to achieve. PIML-OFEM is able to maintain stable iteration under a smaller filter radius, indicating that its displacement field and elemental strain energy prediction errors do not create significant local oscillations or cumulative bias in sensitivity computation, thereby helping to unlock the high-resolution design space.

Further, Fig. 19 compares the results of different mechanical analysis methods used for cantilever beam topology optimization. The geometry, boundary conditions, and load setting of this problem are the same as in Fig. 13, the design domain is discretized into 3061×1531 fine-scale elements, the volume fraction constraint is $v = 0.5$, and the objective remains structural compliance minimization. If classical FEM is directly used for each step of equilibrium analysis, the average single-step analysis time is about 105 s. Fig. 19(a) shows the optimization result obtained based on the traditional substructure method with the boundary linear displacement assumption. Since this method has limited expressive capability for complex boundary deformation and local material heterogeneity, it tends to introduce local analysis errors at substructure interfaces, load application regions, and near material boundaries. To ensure stable optimization, the filter radius needs to be increased to about 5 times the fine element size, with an average single-step analysis time of about 25s.

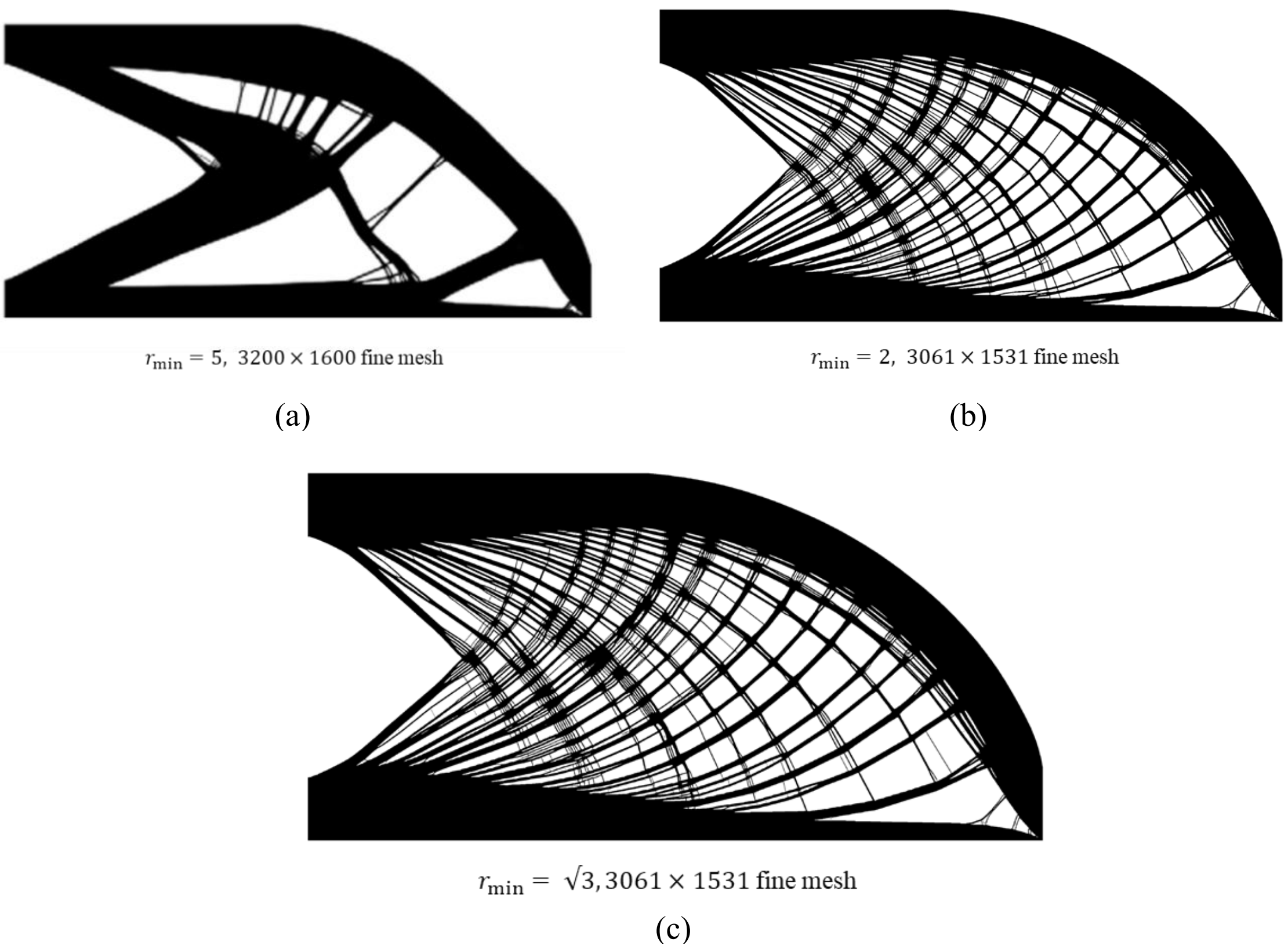

Fig.19. Cantilever-beam topology optimization with different structural-analysis models. (a) Conventional substructure method with linear interpolation of the substructure-boundary displacements. (b) PIML-OFEM with $r_{\min}$ equal to two fine-scale element widths. (c) PIML-OFEM with $r_{\min}$ equal to $\sqrt{3}$ fine-scale element widths.

Fig. 19(b-c) respectively give the optimization results obtained with PIML-OFEM at filter radii of 2 times and $\sqrt{3}$ times the fine element size. The average single-step analysis time of PIML-OFEM is about 33 s, significantly lower than classical FEM, while also being able to obtain stable optimized configurations with richer geometric detail at a smaller filter radius. As the filter radius decreases, fine-scale strut members and local hierarchical features gradually emerge in the structure, indicating

that the proposed method can well preserve the design freedom advantage brought by the high-resolution background mesh. Compared with the traditional boundary-linear-assumption substructure method, although PIML-OFEM has a slightly higher single-step time, it reduces the need for strong filtering introduced by local analysis errors, thereby obtaining higher-quality, higher-detail-resolution optimization results.

In summary, the numerical examples in this section show that PIML-OFEM exhibits good accuracy and efficiency across different material distributions, complex topologies, and topology optimization iteration scenarios. This method alleviates the expressive limitations brought by the linear boundary-displacement assumption of the target substructure through oversampled numerical basis functions, maintains global displacement compatibility through overlapping partition of unity, and reduces the cost of constructing local shape functions through a problem-independent U-Net prediction model. Compared with classical FEM, PIML-OFEM can significantly reduce online linear system solving cost; compared with traditional substructure reduction methods, it can effectively recover local deformation and strain energy distribution while maintaining a low number of degrees of freedom. These results together verify the potential of PIML-OFEM as a computational kernel for large-scale heterogeneous structure analysis and high-resolution topology optimization.

# 6 Conclusion

This paper addresses the difficulty of simultaneously achieving efficiency and accuracy in large-scale heterogeneous structural analysis and high-resolution topology optimization, and proposes an overlapping finite element substructure analysis framework based on problem-independent machine learning, namely PIML-OFEM. The core of this framework is not to replace the overall finite element solution end-to-end with a neural network, but rather, within a low-degree-of-freedom substructure condensation framework, to couple oversampled numerical basis functions, an overlapping finite element compatibility mechanism based on partition of unity, and U-Net local mapping prediction. In doing so, the method retains only the corner node degrees of freedom of the substructure while avoiding the need to directly impose a linear displacement pattern on the substructure's target boundary, providing a unified computational foundation for fine-scale displacement recovery, strain energy computation, and topology optimization sensitivity evaluation under complex heterogeneous materials.

In terms of local reduced-order space construction, this paper uses oversampling techniques to construct numerical basis functions within an expanded region covering the target substructure, and then restricts them to the target substructure. Unlike methods based on linear or higher-order interpolation assumptions on the substructure boundary, this treatment allows the displacement pattern at the target substructure boundary to be naturally induced by the fine-scale material distribution and local mechanical response within the expanded region, thereby improving the low-degree-of-freedom substructure space's ability to express complex local deformation. Through rigid-body mode reproduction conditions and corner node Kronecker-δ normalization, the resulting numerical basis functions maintain consistency of corner node degree-of-freedom interpolation and avoid spurious strain energy under rigid-body motion, providing a stable approximation space for constructing the local condensed stiffness matrix.

In terms of global compatibility, this paper addresses the issue that independently constructed substructure shape functions may produce displacement jumps on shared edges by introducing an overlapping finite element approach based on partition of unity. By designing a specific non-uniform internal substructure mesh, the overlapping of adjacent substructures still forms a regular global Cartesian fine mesh; weight functions are further used to combine local interpolations within the overlapping region in a weighted manner, constructing a globally compatible displacement field. This mechanism allows the oversampled local mapping obtained in Section 2 to enter the global compatible space, thereby improving the reliability of derived quantities such as strain, strain energy, and topology optimization sensitivity. At the same time, the classification into 9 substructure types and 16 overlapping fine element types allows the related integral terms to be precomputed offline and reused during the online stage, so the overlapping compatibility mechanism does not significantly increase runtime assembly complexity.

The numerical examples verify the effectiveness of the proposed method from three perspectives: displacement accuracy, derived response quantities, and topology optimization stability. The cantilever beam and complex hole-structure examples show that PIML-OFEM can achieve displacement responses close to those of classical fine-scale FEM while maintaining a low number of solution degrees of freedom, and can well recover local deformation near hole boundaries, load application regions, and complex material distributions. The elemental strain energy results further show that the oversampled shape functions, the overlapping compatibility mechanism, and a small amount of residual correction can jointly improve the post-processing accuracy of displacement derivative quantities. In the topology optimization example, after embedding PIML-OFEM into the classical SIMP workflow, the optimization process can proceed stably at million-scale fine-element resolution, obtaining optimized configurations with rich fine-scale topological features under a smaller filter radius.

Overall, the advantages of PIML-OFEM are mainly reflected in three aspects. First, without increasing the substructure's boundary control degrees of freedom, it alleviates the limitation on the expressive capability for complex local deformation imposed by the linear assumption on the target substructure boundary. Second, the overlapping compatibility strategy based on partition of unity organizes the independently constructed local oversampled mappings into a globally compatible displacement field, enabling the condensed analysis results to serve tasks such as strain, strain energy, and topology optimization sensitivity that are more sensitive to displacement continuity. Third, the U-Net prediction model significantly reduces the online construction cost of the oversampled numerical basis function, enabling this framework to be used for large-scale heterogeneous structural analysis and high-resolution topology optimization. The regular Cartesian mesh, limited number of element types, and batch neural network inference also facilitate further vectorized and parallel implementation.

There remain several issues in this method worth further investigation. The current validation focuses mainly on two-dimensional linear elastic, small-deformation problems; for three-dimensional structures, geometric nonlinearity, material nonlinearity, and dynamics problems, the construction of oversampled shape functions, integration over overlapping regions, and the output size of the machine learning model would all increase significantly, requiring the development of corresponding efficient implementation strategies. This paper also uses a fixed substructure size and a fixed number of oversampling layers for training and validation; future work could investigate a unified network representation for multi-scale, multi-size substructures, and combine this with error indicators or uncertainty estimation to adaptively trigger local precise oversampling or finite element re-analysis in high-error regions. In addition, combining PIML-OFEM with manufacturing constraints, microstructure optimization, and multi-scale design methods is also expected to further expand its application in rapid analysis and high-performance design of complex engineering structures.

## Conflict of interest

On behalf of all authors, the corresponding author states that there is no conflict of interest.

## Author contributions

**Yilin Guo**: Methodology, Investigation, Formal analysis, Validation, Software, Data curation, Visualization, Funding acquisition, Writing – original draft, Writing – review & editing. **Chang Liu**: Conceptualization, Methodology, Investigation, Formal analysis, Funding acquisition, Writing – original draft, Writing – review & editing. **Zongliang Du**: Formal analysis, Validation, Writing – review & editing. **Jin Liu**: Software. **Jingyu Feng**: Software. **Xinyang Zhang**: Software. **Yang Li**: Software. **Tianxing Yang**: Software. **Changyu Shen**: Resources, Supervision, Writing – review & editing. **Xu Guo**: Conceptualization, Methodology, Project administration, Resources, Funding acquisition, Supervision, Writing – review & editing.

## Acknowledgments

This work was supported by the Science and Technology Key Project of Liaoning Province (2024JH1/11700045), the National Natural Science Foundation of China (124B2042, 12472344), and the Dalian Science and Technology Talent Innovation Support Plan (2024RG001).

# PIML-OFEM：一种基于问题无关机器学习与重叠有限元技术的新型大规模结构分析方法

郭一麟，刘畅，杜宗亮，刘进，冯靖宇，张欣阳，李阳，杨添行，
申长雨，郭旭

**摘要**

高分辨率结构分析与大规模异质结构快速设计，既需要精确的降阶模型，也需要高效的在线计算。在传统多尺度有限元框架下，针对不同材料分布在线构造多尺度形函数面临高昂的计算代价；而前期发展的基于子结构的问题无关机器学习（PIML）技术，其分析精度常受限于子结构边界上假设的位移分布形式。为应对上述挑战，我们提出 PIML-OFEM 一种由问题无关机器学习加速的重叠有限元方法。该方法仅保留每个子结构的角节点自由度，并通过在扩展域上求解局部边值问题，构造与这些自由度相关联的数值基函数。这种降阶模型的构造方式无需在目标子结构边界预设位移假设，显著提升了数值基函数再现子结构局部变形的能力。进一步引入单位分解重叠有限元格式，将子结构上独立构造的数值基函数插值为具有全局连续性的整体基函数。为避免在线求解子结构上的数值基函数，我们训练了一个 UNet，用于学习从子结构杨氏模量分布到其数值基函数的确定性映射。由于该学习过程不依赖于特定的荷载工况与边界条件，所习得的数值形函数可适用于任意同类结构分析与优化问题的求解。数值算例表明，PIML OFEM 获得的位移与单元应变能与细尺度有限元结果高度吻合，其在线计算成本显著低于直接有限元分析，同时计算精度大幅优于基于边界线性插值的 PIML 方法。当将该方法嵌入拓扑优化流程时，PIML-OFEM 能够实现小滤波半径的高分辨率优化，并保留包括类秩 2 微结构局部模式在内的细尺度结构特征。通过融合基于力学原理的局部独立降阶、全局连续耦合以及问题无关的机器学习技术，PIML-OFEM 为大尺度异质结构分析与高分辨率拓扑优化建立了一种高效的人工智能赋能数值计算新范式。